\documentclass{amsart}

\RequirePackage{color}

\def\smallskip{\vskip\smallskipamount}
\def\medskip{\vskip\medskipamount}
\def\bigskip{\vskip\bigskipamount}

\usepackage{amsmath,amsthm,bm,mathrsfs,amscd,tikz}
\usepackage{ytableau}
\usepackage{comment}
\usepackage{mathdots}
\usepackage[all]{xy}
\usepackage{graphicx}
\usepackage[colorinlistoftodos]{todonotes}
\usepackage{enumerate}
\usepackage{feynmp-auto}
\usepackage[english]{babel}
\usepackage[utf8]{inputenc}
\usepackage{float}
\usepackage{tikz}
\usepackage{tikz-cd}
\usepackage{amssymb}
\usepackage{mathtools}
\usepackage{soul}
\usetikzlibrary{decorations.pathmorphing}

\usepackage[utf8]{inputenc}

\newtheoremstyle{thmstyle}{}{}{\itshape}{}{\bfseries}{ }{5pt}{}
\newtheoremstyle{exstyle}{}{}{}{}{\bfseries}{ }{5pt}{}
\newtheoremstyle{defstyle}{}{}{}{}{\bfseries}{ }{5pt}{}
\newtheoremstyle{remstyle}{}{}{}{}{\bfseries}{ }{5pt}{}

\theoremstyle{thmstyle}
\newtheorem{thm}{Theorem}[section]
\newtheorem{theorem}[thm]{Theorem}

\newtheorem{lemma}[thm]{Lemma}

\newtheorem{proposition}[thm]{Proposition}

\newtheorem{corollary}[thm]{Corollary}

\newtheorem{conjecture}{Conjecture}[section]

\theoremstyle{exstyle}

\newtheorem{example}[thm]{Example}

\theoremstyle{defstyle}

\newtheorem{definition}[thm]{Definition}

\newtheorem{def-prop}[thm]{Definition-Proposition}
\newtheorem{def-lem}[thm]{Definition-Lemma}
\newtheorem{rem-convention}[thm]{Remark-Convention}
\newtheorem{def-note}[thm]{Definition-Notation}

\theoremstyle{remstyle}

\newtheorem{remark}[thm]{Remark}

\theoremstyle{remstyle}

\newcommand{\Hom}{\operatorname{Hom}}
\newcommand{\Fac}{\operatorname{Fac}}
\newcommand{\Ext}{\operatorname{Ext}}

\newcommand{\taurigid}{\operatorname{\tau-rigid}}

\DeclareMathOperator*{\rad}{rad}

\DeclareMathOperator*{\modu}{mod}

\DeclareMathOperator*{\proj}{proj}

\newcommand{\Z}{\mathcal{Z}}

\DeclareMathOperator*{\ind}{ind}

\DeclareMathOperator*{\tors}{tors}

\DeclareMathOperator*{\End}{End}

\DeclareMathOperator*{\Filt}{filt}

\DeclareMathOperator*{\brick}{brick}

\DeclareMathOperator*{\node}{node}

\DeclareMathOperator*{\pd}{pd}

\DeclareMathOperator*{\GL}{GL}

\DeclareMathOperator*{\rep}{rep}

\DeclareMathOperator*{\Irr}{Irr}

\DeclareMathOperator*{\coker}{coker}
\DeclareMathOperator*{\Ima}{Im}

\newcommand{\cupdot}{\mathbin{\mathaccent\cdot\cup}}

\makeatletter
\newcommand{\doublewidetilde}[1]{{%
  \mathpalette\double@widetilde{#1}%
}}
\newcommand{\double@widetilde}[2]{%
  \sbox\z@{$\m@th#1\widetilde{#2}$}%
  \ht\z@=.9\ht\z@
  \widetilde{\box\z@}%
}
\makeatother

\subjclass [2020]{16P10,16G20,16D80,16G60}

\begin{document}

\title{Geometric interactions between bricks and $\tau$-rigidity}
\author[Kaveh Mousavand, Charles Paquette]{Kaveh Mousavand, Charles Paquette} 
\address{Unit of Representation Theory and Algebraic Combinatorics, Okinawa Institute of Science and Technology, Okinawa, Japan}
\email{mousavand.kaveh@gmail.com}
\address{Department of Mathematics and Computer Science, Royal Military College of Canada, Kingston ON, Canada}
\email{charles.paquette.math@gmail.com}

\begin{abstract}
For finite-dimensional algebras over algebraically closed fields, we consider two fundamental classes of modules and their geometric counterparts: bricks and $\tau$-rigid modules, as well as brick components and $\tau$-regular components. We then apply our results in the study of some open conjectures. 

First, we investigate the situation where every brick is $\tau$-rigid. We prove that this occurs exactly when the algebra is \emph{locally representation-directed}; a family of algebras introduced by Dräxler in the 1990s, which are always representation-finite. 
Then, we adopt a geometric perspective and analyze the brick and $\tau$-regular components of module varieties. In this greater generality, we establish new properties of such components. Inspired by some recent conjectures, we apply our results to the study of minimal brick-infinite algebras. Along the way, we construct some limits of rigid $g$-vectors, under a condition that we call the \emph{$\tau$-convergence property}. This construction is novel and, in certain cases, yields an integral $g$-vector lying outside the $\tau$-tilting fan (a.k.a. $g$-vector fan). We show how our results provide new tools to the study of some open conjectures and particularly illustrate that for $E$-tame algebras.

\end{abstract}

\maketitle

\tableofcontents

\section{Introduction}

To motivate our problems, put them into a larger perspective, and summarize our main results, we briefly describe the setting and recall some background. The notation and terminology which are not explicitly introduced here are recalled in the following sections. For rudimentary materials in representation theory of algebras, we refer to \cite{ASS}. 
Henceforth, $A$ denotes a finite dimensional associative unital algebra of rank $n$ over an algebraically closed field $k$. Without loss of generality, $A$ is assumed to be basic and connected, and hence $A$ can be presented by a connected bound quiver $(Q,I)$; that is, $A\simeq kQ/I$ is an isomorphism of $k$-algebras, and every $A$-module can be seen as a representation of $(Q,I)$. 
Let $\modu A$ denote the category of finitely generated left $A$-modules. Unless specified otherwise, a module always means an object in $\modu A$, considered up to isomorphism. By $\ind(A)$ denote the set of all (isomorphism classes) of indecomposable objects in $\modu A$. Recall that $A$ is said to be \emph{representation-finite} if $\ind(A)$ is a finite set. 

\medskip

Many fundamental studies in classical representation theory were motivated by the classification, size, and distribution of indecomposable modules. In that regard, a better understanding of the indecomposable modules over representation-infinite algebras was among the ultimate goals. 
Some of these phenomena were formulated in terms of Brauer-Thrall Conjectures (now theorems) formally posed in the late 1950's, which inspired several directions of research in the second half of the $20^{th}$ century (for details, see \cite{Bo} and references therein). 
Meanwhile, for a given algebra $A$, it is known that the rich dictionary between the algebraic and geometric aspects of the representation theory of $A$ provides a range of powerful tools in the study of many classical and modern problems. For instance, the Second Brauer-Thrall Theorem combined with Voigt's Lemma imply that ``\emph{If (almost) all modules in $\ind(A)$ are rigid, then $A$ is representation-finite} (Namely, if for all, but at most finitely many, $X  \in \ind(A)$ we have $\Ext^1_A(X,X)=0$, then $|\ind(A)|<\infty$.)". 

\medskip

Recently, new notions of rigidity and various interactions between different sets of indecomposables have received a lot of attention. For example, thanks to $\tau$-tilting theory introduced in \cite{AIR}, the connections between indecomposable $\tau$-rigid modules and other subsets of $\ind(A)$ are studied extensively. 
Recall that $M$ in $\modu A$ is \emph{$\tau$-rigid} if $\Hom_A(M,\tau M)=0$, where $\tau M$ denotes the Auslander-Reiten translate of $M$ in $\modu A$. By the Auslander-Reiten duality, every $\tau$-rigid module is rigid (but the converse is not necessarily true). Let $\textit{i}\taurigid(A)$ denote the set of all $\tau$-rigid modules belonging to $\ind(A)$. Also recall that $X$ in $\modu A$ is a \emph{brick} if $\End_A(X)$ is a division algebra. By $\brick(A)$ we denote the set of all (isomorphism classes) of bricks in $\modu A$. In our setting, $\brick(A)$ consists of exactly those modules $X$ in $\modu A$ for which $\End_A(X)\simeq k$. That being the case, bricks are also known as Schur representation. For a recent report on some studies of bricks, see \cite{MP4}.

\medskip

In his PhD dissertation, the first-named author posed two conjectures, the stronger of which can be viewed as the brick-analogue of the Second Brauer-Thrall Conjecture (see Conjecture \ref{2nd bBT conj.}). Meanwhile, it is notable that this conjecture also immediately implies a brick-analogue of the above-mentioned well-known results on the rigidity of almost all indecomposable modules. More precisely, phrased in a less technical language, \cite[Conjecture 6.0.1]{Mo1} implies ``\emph{If (almost) all modules in $\brick(A)$ are rigid, then $A$ is brick-finite} (Namely, if for all, but at most finitely many, $X \in \brick(A)$ we have $\Ext^1_A(X,X)=0$, then $|\brick(A)|<\infty$.)". Even this simpler implication is unknown in full generality and remains an open problem!

\medskip

Inspired by the above problem and some other conjectures recalled below, in this work we primarily focus on the stronger notion of $\tau$-rigidity and study the algebraic and geometric properties of bricks and some related indecomposable modules and irreducible components. 
In the earlier part of this paper, we make a careful analysis of the interactions between the three sets of indecomposables in $\modu A$. Our first theorem shows that a ``two-out-of-three" phenomenon holds for the sets $\ind(A)$, $\brick(A)$, and $\textit{i}\taurigid(A)$. 
Before stating our first main result, we recall that in \cite{DIJ} the authors established an injective map from $\textit{i}\taurigid(A)$ to $\brick(A)$, known as the ``brick-$\tau$-rigid correspondence" (see Theorem \ref{Thm of DIJ}). 
In Section \ref{Section:new char. of loc-rep-directed}, we show that $\brick(A) \subseteq \textit{i}\taurigid(A)$ holds if and only if $A$ is ``locally representation-directed". First introduced in \cite{Dra}, this family generalizes representation-directed algebras, and our first theorem gives new characterizations of these interesting algebras.

\begin{theorem}[Theorem \ref{Thm:new characterization of locally rep-directed} and Corollary \ref{Cor: locally rep-directed}]\label{Thm 1.1}

For an algebra $A$, if two of the three sets $\ind(A)$, $\brick(A)$, and $\textit{i}\taurigid(A)$ are equal, then all the three sets coincide. 
Consequently, the following are equivalent: 
\begin{enumerate}
    \item $A$ is locally representation-directed;
    \item Every brick is $\tau$-rigid;
     \item $\ind(A)=\brick(A)=\textit{i}\taurigid(A)$;
    \item The brick-$\tau$-rigid correspondence is the identity map;
\end{enumerate}
\end{theorem}

Note that the assumption on the $\tau$-rigidity of bricks in the above theorem cannot be replaced by the weaker condition of rigidity (see Example \ref{Examples}). 
Meanwhile, our results immediately imply that if every brick in $\modu A$ is $\tau$-rigid, $A$ is representation-finite, thus brick-finite. This particularly highlights the connection to the open problem noted above (i.e., whether the rigidity of all bricks implies brick-finiteness). 

\medskip

Now we turn our attention to two challenging open conjectures that closely relate to the scope of our algebro-geometric treatment of the representation varieties. 
The first conjecture below is nowadays often known as Second brick-Brauer-Thrall Conjecture, abbreviated as 2nd bBT Conjecture, and concerns the behavior of those algebras which are brick-infinite, that is, when $\brick(A)$ is an infinite set.
This conjecture posed in a more technical formulation first appeared in 2019, in the preprint of \cite{Mo2}, and also as \cite[Conjecture 6.0.1, part (2)]{Mo1}. 

\medskip

\begin{conjecture}[2nd bBT Conj.] \label{2nd bBT conj.}
If $A$ is brick-infinite, then for some integer $d$, there are infinitely many (non-isomorphic) bricks of dimension $d$.    
\end{conjecture}

To understand the next conjecture, one needs the notion of $E$-finiteness, which can be expressed in different ways. One formulation is that $A$ is an \emph{$E$-finite} algebra provided its $\tau$-tilting fan is complete in the rational Grothendieck group $K_0(\proj A)_\mathbb{Q}$. For the definition of the $\tau$-tilting fan and more details on $E$-finiteness, see Subsection \ref{section tau-tilting fan}. It is known that a brick-finite algebra is necessarily $E$-finite (a consequence of Theorem \ref{Thm: complete tau-fans}). The conjecture below is about the converse, first posed in 2017, by Demonet \cite[Question 3.49]{De}.

\begin{conjecture}[Demonet's Conj.] \label{Demonet Conj intro} 
If $A$ is brick-infinite, then $A$ is not $E$-finite.
\end{conjecture} 

Although the above conjectures are well-known for hereditary algebras and some other families, we remark that both conjectures are still open in full generality. Meanwhile, there has been some important progress in the past few years that settled them in some cases and developed substantial reductions (for details, see \cite{Mo1, Mo2, MP1, MP2, MP3, MP4, STV, Pf1, Pf2} and references therein). 

\medskip

Inspired by our results in Theorem \ref{Thm 1.1}, and in connection with Conjectures \ref{2nd bBT conj.} and \ref{Demonet Conj intro}, in the second part of our paper we consider the geometric analogues of the three sets  $\ind(A)$, $\brick(A)$, and $\textit{i}\taurigid(A)$. 
More precisely, for algebra $A$, we consider the representation varieties of $A$ and their irreducible components. Let $\Irr(A)$ denote the set of all irreducible components in all representation varieties $\rep(A,\textbf{d})$, where $\textbf{d}$ is an arbitrary dimension vector of $A$ (see Section \ref{Preliminaries}).
As shown in \cite{CBS}, the study of components in $\Irr(A)$ heavily relies on the knowledge of indecomposable components, namely, those $\Z$ in $\Irr(A)$ which contain a non-empty open subset consisting of indecomposable modules. Let $\Irr^i(A)$ denote the set of all indecomposable components in $\Irr(A)$. Moreover, by $\Irr^b(A)$ we denote the set of brick components in $\Irr(A)$, which are those irreducible components of $A$ that contain a brick (it then follows that the bricks form a non-empty open set in the component). Finally, by $\Irr^{i\tau}(A)$ we denote the set indecomposable components which are $\tau$-regular (for definition, see Section \ref{Preliminaries}).

\medskip

We view the sets $\Irr^{i}(A)$, $\Irr^{b}(A)$ and $\Irr^{i\tau}(A)$ as the geometric counterparts of the three algebraic sets $\ind(A)$, $\brick(A)$, and $\textit{i}\taurigid(A)$, treated in Theorem \ref{Thm 1.1}. However, this analogy is not just a verbatim translation of the algebraic setting to the geometric one.
First, note that $A$ may admit irreducible components that contain indecomposables but they are not indecomposable component. In particular, not every module of $\ind(A)$ gives rise to a component in $\Irr^i(A)$.
Secondly, a brick component may simultaneously contain several parameters of non-isomorphic bricks. In fact, for each $\Z$ in $\Irr^{b}(A)$, either $\Z$ is the closure of a single orbit of a brick, or else $\Z$ contains infinitely many orbits of bricks.  Thus, we have a surjective map from $\brick(A)$ to $\Irr^b(A)$.
Lastly, observe that the modules in $\textit{i}\taurigid(A)$ are in bijection with those components $\Z$ in $\Irr^{i\tau}(A)$ which are orbit closures. Hence, there is an injective map from $\textit{i}\taurigid(A)$ to $\Irr^{i\tau}(A)$. In fact, $A$ is $E$-finite if and only if this map is bijective. Thus, for an algebra $A$, there may exist components in $\Irr^{i\tau}(A)$ which are not an orbit closure of a $\tau$-rigid module (see Section \ref{Section:gen tau-reduced components}).

\medskip

Before we state the next result, note that Theorem \ref{Thm 1.1} shows that the equality $\brick(A)=\textit{i}\taurigid(A)$ captures some important properties of the algebras under consideration. Our next theorem indicates how the equality $\Irr^{b}(A)=\Irr^{i\tau}(A)$ can be seen as a geometric generalization of the locally representation-directed algebras. This is further exploited in Section \ref{Section:geometric generalization of loc-rep-directed}. Recall $A$ is said to be \emph{minimal brick-infinite} if $A$ is brick-infinite but every proper quotient algebra of $A$ is brick-finite. 
For the decisive role of minimal brick-infinite algebras in the study of several open problems, see \cite{MP1,MP3,MP4}. 
We remark that, to verify Conjectures \ref{2nd bBT conj.} and \ref{Demonet Conj intro} in full generality, it is sufficient to settle them for minimal brick-infinite algebras. 

\begin{theorem}\label{thm: Geometric two-out-of-three}

Let $A$ be an algebra for which two of the three sets $\Irr^i(A)$, $\Irr^b(A)$, and $\Irr^{i\tau}(A)$ are equal. Then, $\textit{i}\taurigid(A) \subseteq \brick(A)$, and furthermore we have:
\begin{enumerate}
    \item $A$ is brick-finite if and only if $A$ is locally representation-directed.
    \item If $A$ is minimal brick-infinite then $A$ admits (at least) one infinite family of bricks in the same dimension; that is, Conjecture \ref{2nd bBT conj.} holds for $A$.
\end{enumerate}
\end{theorem}

The ``two-out-of-three'' assumption in Theorems \ref{Thm 1.1} and \ref{thm: Geometric two-out-of-three} may initially seem to hold only for a very small family of algebras. In the algebraic setting of Theorem \ref{Thm 1.1}, this indeed classifies a subfamily of representation-finite algebras, but in the geometric setting of Theorem \ref{thm: Geometric two-out-of-three}, we have a much larger family. In fact, in Section \ref{Section:geometric generalization of loc-rep-directed} , we show that for any $n>1$, there are infinitely many non-isomorphic algebras of rank $n$ for which we have the equalities $\Irr^i(A)=\Irr^b(A)=\Irr^{i\tau}(A)$, and $A$ can be chosen to be representation-finite, or tame type, as well as of wild type.

\medskip

To motivate our next result, note that although Conjecture \ref{Demonet Conj intro} has proved to be important in various areas of modern representation theory and has received attention from different researchers, to our knowledge, there is no general tool for constructing a rational $g$-vector which falls outside of the $\tau$-tilting fan of a given brick-infinite algebra. 
Equivalently, starting from a brick-infinite algebra $A$, there is no systematic approach which, under certain conditions, allows one to construct an irreducible $\tau$-regular component without an open orbit. 
In Section \ref{Section:min-brick-inf}, we introduce such a tool which leads to the desired results under a condition that we call the ``$\tau$-convergence property". Over minimal brick-infinite algebras, this construction gives some new results on this open problem. 
Observe that, as discussed in Remark \ref{Rem: Reduction of Demonet's Conj}, to settle Conjecture \ref{Demonet Conj intro}, it suffices to verify it over minimal brick-infinite algebras.

\begin{theorem} [Theorem \ref{thm:rational_g_vector}]
Every minimal brick-infinite algebra $A$ having the $\tau$-convergence property is $E$-infinite; that is, Conjecture \ref{Demonet Conj intro} holds for $A$.

\end{theorem}

To prove the above result, we develop some techniques that combine some important homological properties of minimal brick-infinite algebras with properties of a quadratic form on $g$-vectors defined via the Cartan matrix of $A$. We suspect that the $\tau$-convergence property always holds for $g$-tame and minimal brick-infinite algebras. 
In fact, we obtain the following interesting corollary for tame algebras, showing that our $\tau$-convergence property provides a unified approach to the study of both open Conjectures \ref{2nd bBT conj.} and \ref{Demonet Conj intro}. 
We note that the final assertion of the next corollary follows from a reduction theorem \cite[Theorem 1.4]{MP1} in our earlier work.

\begin{corollary}[Remark \ref{2nd bBT and Demonet's conj for tame algebras} and Prop. \ref{Prop: 2nd bBT and  tau-convergence for tame algebras}]
For tame algebras, the following are equivalent:

\begin{enumerate}
    \item Conjecture \ref{2nd bBT conj.} holds;
    \item Conjecture \ref{Demonet Conj intro} holds;
    \item Every minimal brick-infinite (tame) algebra has the $\tau$-convergence property.
\end{enumerate}
Moreover, to verify one (therefore, all) of the above conditions, it suffices to consider those minimal brick-infinite tame algebras over which almost all bricks are faithful.
\end{corollary}

Due to the central role of minimal brick-infinite algebras in our studies, in Section \ref{Section:min-brick-inf} we analyze their $\tau$-tilting fan and wall-and-chamber structure. 
Before we summarize these results, recall that to each $g$-vector $\theta$, a unique $\tau$-regular component can be associated in $\Irr(A)$, which we denote by $\Z(\theta)$ (see Theorem \ref{Plamondon's Thm}). 
Let us note that $\theta$ belongs to the $\tau$-tilting fan $\mathcal{F}_A$ if and only if $c(\Z(\theta))=0$. 
We say that an irreducible component $\Z$ in $\Irr(A)$ is \emph{faithful} if, for each nonzero ideal $J$ in $A$, there exists $Z \in \Z$ such that $JZ\neq 0$. For each irreducible component $\Z$, we observe that the set of faithful modules forms an open subset (see Prop. \ref{prop: faithful is open}). For a more refined version of the following result, see Proposition \ref{Prop:measure zero}. 

\begin{proposition}
Let $A$ be a minimal brick-infinite algebra. Then, the $g$-vectors $\theta$ for which $\Z(\theta)$ is not faithful is a finite union of some chambers and walls.
\end{proposition}

In the final section, we consider the torsion classes induced by the components in $\Irr(A)$ and connect their behavior to the study of Conjectures \ref{2nd bBT conj.} and \ref{Demonet Conj intro}. We then apply our results to $E$-tame algebras and derive new results on some open problems. 
We particularly make an important observation on those $E$-tame algebras over which almost all bricks are $\tau$-rigid (compare with Theorem \ref{Thm 1.1}). To put this into a better perspective, note that there exist algebras for which almost all bricks are $\tau$-rigid, but the algebra is representation-infinite. Hence, it is natural to study the important consequence of Conjecture \ref{2nd bBT conj.} noted before, that is, ``\emph{If almost all modules in $\brick(A)$ are $\tau$-rigid, then $A$ is brick-finite}". In the absence of a full proof for this open problem, in Proposition \ref{prop: E-tame almost all brick tau-rigid} we verify that for all $E$-tame algebras, hence for all tame algebras. 
Our results are summarized in the next theorem; we note that part $(1)$ has been independently shown in \cite{Pf2}, using different techniques.

\begin{theorem}[Corollary \ref{Cor: E-tame and bricks} and Prop. \ref{prop: E-tame almost all brick tau-rigid}]
Let $A$ be an $E$-tame algebra. Then,
\begin{enumerate}
    
    \item If $A$ is not $E$-finite, it has an infinite family of bricks of the same dimension. Consequently, for $E$-tame algebras, Conjecture \ref{Demonet Conj intro} implies Conjecture \ref{2nd bBT conj.}.

    \item If almost all bricks of $A$ are $\tau$-rigid, then $A$ is brick-finite. 
\end{enumerate}
\end{theorem}

\section{Preliminaries, setting and notations}\label{Preliminaries}
We always work with finite dimensional algebras over an algebraically closed field $k$. In particular, $A$ denotes a finite dimensional associative $k$-algebra with multiplicative identity, and without loss of generality, $A$ is assumed to be connected, basic, and of rank $n$. Hence, we can always consider a fixed a presentation of $A$ in terms of a bound quiver $(Q,I)$, where $Q$ is a quiver uniquely associated to $A$, and $I$ is an admissible ideal in the path algebra $kQ$ such that $A\simeq kQ/I$ is an isomorphism of $k$-algebras. Consequently, as we do throughout the text, modules over $A$ can be viewed as representations of the bound quiver $(Q,I)$, and vice versa.
By $\modu A$, we denote the category of finitely generated left $A$-modules. Unless specified otherwise, a module always mean an object in $\modu A$, considered up to isomorphism. Let $\ind(A)$ denote the set of all (isoclasses of) indecomposable modules in $\modu A$, and recall that $A$ is said to be \emph{representation-finite} if $\ind(A)$ is a finite set.

\medskip
For $M$ in $\modu A$, let $|M|$ denote the number of non-isomorphic (and non-zero) indecomposable summands of $M$ in its Krull-Schmidt decomposition. In particular, $|A|=n$.
For $M$, the corresponding representation of the bound quiver $(Q,I)$ is formally denoted by a pair $(\{M_i\}_{i\in Q_0}, \{M(\alpha)\}_{\alpha \in Q_1})$, where $Q_0 = \{1, 2, \ldots, n\}$ and $Q_1$ are respectively the vertex set and arrow set of $Q$, each $M_i$ is a finite dimensional $k$-vector space, and for each arrow $\alpha:i\rightarrow j$ in $Q_1$, we have $M(\alpha): M_i \rightarrow M_j$. Let $\textbf{d}_M$ denote the dimension vector of $M$, which is the vector $(\dim_k(M_i))_{i \in Q_0}$.
For $i \in Q_0$, by $S_i$ we denote the simple $A$-module associated to $i$, and we let $P_i$ and $I_i$ respectively denote the projective and injective indecomposable $A$-modules corresponding to $i \in Q_0$. In particular, observe that $P_i = Ae_i$ is the projective cover of $S_i$, and $S_i = Ae_i/{\rm rad}(A)e_i$, where ${\rm rad}(A)$ denotes the Jacobson radical of $A$.
Thus, $A= Ae_1 \oplus \cdots \oplus Ae_n$, where $\{e_1, \dots, e_n\}$ is a complete set of primitive orthogonal idempotents in $A$.
All the standard notations and rudimentary materials in representation theory used below can be found in \cite{ASS}.

\medskip
Let $\proj A$ denote the full exact subcategory of $\modu A$ whose objects are projective modules and by $K^b({\rm proj}A)$ denote the homotoy category of bounded complexes of $\proj A$. 
We use the natural identification of the Grothendieck groups $K_0(K^b({\rm proj}A))$ and $K_0(\rm proj A)$ and consider the basis $\{[P_1], \ldots, [P_n]\}$ represented by the indecomposable projective modules. By the isomorphism $K_0(\rm projA) \cong \mathbb{Z}^n$, the element $[P_i]$ corresponds to the canonical basis vector of $\mathbb{Z}^n$ which is $1$ at the $i$-th entry and zero elsewhere. We use $K_0(\rm proj A)_{\mathbb{Q}}$ and $K_0(\rm proj A)_{\mathbb{R}}$ to respectively denote $K_0(\rm proj A) \otimes_{\mathbb{Z}} \mathbb{Q}$ and $K_0(\rm proj A) \otimes_{\mathbb{Z}} \mathbb{R}$. 
On the other hand, we have $K_0(D^b(\modu A))\simeq \mathbb {Z}^n$, where $D^b(\modu A)$ denotes the bounded derived category of $\modu A$. Here, we have the standard basis $\{[S_1], \ldots, [S_n]\}$, given by the isomorphism classes of simple modules.
For $M$ in $\modu A$, we can identify its dimension vector $\textbf{d}_M$ with the corresponding element $[M]$ in $K_0(D^b(\modu A))$. 

\medskip

Each element of $K_0(\rm proj A)$ is often called a \emph{$g$-vector}. Every $g$-vector, say $v$, can be uniquely written as a difference of two vectors with non-negative entries, that is, $v$ = $v^+ - v^-$, where $v^+ = (v'_i)$ such that $v_i' = {\rm max}\{0, v_i\}$, and $v^-=(v''_i)$ given by $v''_i = -{\rm min}\{0, v_i\}$. Thus, $P(v^+) := \oplus_{i=1}^n P_i^{(v^+)_i}$ and $P(v^-) := \oplus_{i=1}^n P_i^{(v^-)_i}$ are well-defined projective modules. Therefore, to the $g$-vector $v$, we naturally associate the $2$-term complex $P(v^-) \to P(v^+)$ in $K^b({\proj}(A))$. 
Note that, $\Hom_A(P(v^-), P(v^+))$ can be seen as an affine variety, consisting of points $f = (f_i)$ in the vector space $\prod_{i=1}^n\Hom_k(P(v^-)_i, P(v^+)_i)$ where the linear conditions $f_jP(v^-)(\alpha) = P(v^+)(\alpha)f_i$ are simultaneously satisfied, for all $\alpha: i \to j$ in $Q_1$. Alternatively, each $f = (f_i)$ in $\Hom_A(P(v^-), P(v^+))$ can be seen as a matrix of size $|v^+| \times |v^-|$, where for a vector $u$, by $|u|$ we denote the sum of its entries. A given entry of this matrix represents the maps $P \to R$ between some indecomposable projective modules $P$ and $R$. In particular, if $P = P_i$ and $R = P_j$, there is a basis
$B_{ij}=\{p_{ij}^1, \ldots, p_{ij}^{s_{ij}}\}$ of $\Hom_A(P_i, P_j) = e_iAe_j$, and each entry corresponding to $P\to R$ is a linear combinations of the elements in $B_{ij}$.

\subsection{Irreducible components of module varieties}\label{Subsection: Generically tau-reduced components}
For an algebra $A=kQ/I$ of rank $n$, and a fixed dimension vector $\textbf{d}:=(d_i)_{i=1}^n \in \mathbb{Z}^{Q_0}_{\geq 0}$, by $\rep(A,\textbf{d})$ we denote the affine variety parametrizing the representations of $(Q,I)$ of dimension vector $\textbf{d}$. The general linear group $\GL(\textbf{d}):=\GL(d_1)\times \cdots \times\GL(d_n)$ acts on $\rep(A,\textbf{d})$ via conjugation.
For $M \in \rep(A,\textbf{d})$, by $\mathcal{O}_M$ we denote the $\GL(\textbf{d})$-orbit of $M$. More generally, for $X \in \modu A$, let $\mathcal{O}_X$ denote the orbit of $X$ in the representation variety that $X$ belongs to. Through this dictionary, $X$ and $Y$ in $\modu A$ lie in the same orbit if and only if they are isomorphic.
\medskip

Let $\Irr(A,\textbf{d})$ denote the set of all irreducible components of the variety $\rep(A,\textbf{d})$, and  by $\Irr(A)$ denote the union of all $\Irr(A,\textbf{d})$, where $\textbf{d}$ runs through $\mathbb{Z}^{n}_{\geq 0}$. 
Recall that $\mathcal{Z} \in \Irr(A)$ is said to be \emph{indecomposable} if there is a non-empty open subset $\mathcal{U} \subset \Z$ consisting of indecomposable representations.
For $\mathcal{Z} \in \Irr(A)$, it is well-known that the set of all bricks in $\Z$ form an open (possibly empty) subset in $\Z$.
In particular, $\Z$ is called a \emph{brick component} if it contains a brick. Let $\Irr^i(A)$ denote the set of all indecomposable components in $\Irr(A)$, and by $\Irr^b(A)$ denote the set of all brick components. Evidently, $\Irr^b(A)\subseteq \Irr^i(A)$. 

\medskip
For each $\mathcal{Z} \in \Irr(A)$, define $$c(\Z):=\min \{\dim(\Z) - \dim(\mathcal{O}_Z) \,| Z \in \Z \}.$$ 
In other words, $c(\Z)$ denotes the generic number of parameters of $\Z$.
Moreover, let $h(\mathcal{Z})$ denote the minimal value of $\Hom_A(Z, \tau Z)$, where $Z \in \Z$. We observe that $c(\Z)$, as well as $h(\mathcal{Z})$, are generically attained on an open and dense subset of $\Z$. Moreover, in general we always have $c(\mathcal{Z})\leq h(\mathcal{Z})$.
An irreducible component $\Z$ is called \emph{$\tau$-regular} if $c(\Z)= h(\Z)$.
Such components were first considered in \cite{GLS}, and later characterized in \cite{Pl}.
We will see that (Theorem \ref{Plamondon's Thm}) that $\tau$-regular components are in correspondence with some
$g$-vectors. In particular, if $\Z$ is $\tau$-regular, then there is a $g$-vector $v$ corresponding to it. Conversely, if $v$ is any $g$-vector, then there is a unique $\tau$-regular component, denoted $\Z(v)$ corresponding to it.
Similar to our notations above, by $\Irr^{\tau}(A,\textbf{d})$ we denote the set of all $\tau$-regular components in $\rep(A,\textbf{d})$, and $\Irr^{\tau}(A)$ denotes the union of all $\Irr^{\tau}(A,\textbf{d})$, for all $\textbf{d} \in \mathbb{Z}^{n}_{\geq 0}$. Finally, by $\Irr^{i\tau}(A)$ denote the set of all indecomposable components which are $\tau$-regular, that is $\Irr^{i\tau}(A):=\Irr^{i}(A)\cap \Irr^{\tau}(A)$.

\medskip

\subsection{The $\tau$-tilting fan} \label{section tau-tilting fan}
As before, let $A$ of rank $n$, and recall that $M$ in $\modu A$ is called \emph{$\tau$-rigid} if $\Hom(M,\tau M)=0$. Then, a basic \emph{$\tau$-rigid pair} $(M,P)$ consists of a $\tau$-rigid module $M$, and a projective module $P$ in $\modu A$ such that $\Hom_A(P,M)=0$, and both $M$ and $P$ are basic (i.e, in their Krull-Schmidt decomposition, no non-zero indecomposable module appears with multiplicity more than one).
For any such basic $\tau$-rigid pair, the $g$-vectors of the indecomposable summands of $M$ and $P$ form a linearly independent set of elements in $K_0(\rm proj A)_{\mathbb{R}} = \mathbb{R}^n$ and they specify the rays of a polyhedral cone $C(M, P)$ in $\mathbb{R}^n$ associated to $(M,P)$. 
In particular, $C(M, P)$ is a cone of maximal dimension if and only if $(M,P)$ is a support $\tau$-tilting pair, that is, if we have $|M|+|P|=n$. We refer to \cite{AIR} and \cite{DIJ} for more details.

\medskip

The following theorem uses the $\tau$-tilting theory in $\modu A$ to describe the geometry of a fan associated to $A$, which we henceforth call the \emph{$\tau$-tilting fan} of $A$.
In particular, the following result follows from the study of $\tau$-rigid pairs in \cite{AIR}.
\begin{theorem}[\cite{AIR}] \label{Thm: tau-fan}
Let $A$ be an algebra of rank $n$. Then,
the $g$-vectors of indecomposable $\tau$-rigid modules in $\modu A$ are the rays
of a fan $\mathcal{F}_A$ in $\mathbb{R}^n$ whose maximal cones are given by the $g$-vectors of non-isomorphic $\tau$-tilting pairs $(M,P)$. Moreover, $\mathcal{F}_A$ is essential, rational, polyhedral and simplicial.
\end{theorem}

Henceforth, by $\mathcal{F}_A$, or simply $\mathcal{F}$ when there is no risk of confusion,  we denote the $\tau$-tilting fan of $A$ described in the above theorem. We will regard this fan as living in $K_0(\proj(A))_\mathbb{R}$, which is identified with $\mathbb{R}^n$. We can also regard it as a fan in $\mathbb{Q}^n$ but in that case, we will explicitly say so. 
Moreover, $\mathcal{F}^{\circ}$ denotes the union of interior of the maximal cones of $\mathcal{F}$. 

\begin{remark}
The fan $\mathcal{F}$ described above is often seen as the counterpart of the $g$-vector fans of cluster algebras and sometimes it is called the ``$g$-vector fan" of $A$. However, we avoid this term and simply call it the \emph{$\tau$-tilting fan} of $A$,  particularly because $A$ may have $g$-vectors which do not belong to $\mathcal{F}$ (see Theorem \ref{Thm: complete tau-fans}).
\end{remark}

Recall that a fan $\mathcal{G}$ in $\mathbb{R}^n$ is called \emph{complete} if $\mathcal{G}=\mathbb{R}^n$. 

\begin{theorem}[\cite{DIJ}, \cite{As2}]\label{Thm: complete tau-fans}
For an algebra $A$, the $\tau$-tilting fan $\mathcal{F}$ is complete in $\mathbb{R}^n$ if and only if $A$ is $\tau$-tilting finite.
\end{theorem}

We can also think of the $\tau$-tilting fan as living in the rational Grothendieck group $K_0(\proj A)_\mathbb{Q}$, in which case we will use the notation $\mathcal{F}_\mathbb{Q}$. It follows from the previous theorem that when $A$ is brick-finite, then $\mathcal{F}_\mathbb{Q}$ is complete in $K_0(\proj A)_\mathbb{Q}$. We say that $A$ is \emph{$E$-finite} provided $\mathcal{F}_\mathbb{Q}$ is complete in $K_0(\proj A)_\mathbb{Q}$. 

\begin{remark}\label{Remark: on Demonet's conjecture}
From Theorem \ref{Thm: complete tau-fans} and some other known results (see Theorem \ref{Thm of DIJ}), it follows that $A$ is brick-infinite if and only if there are $g$-vectors in the ambient space $\mathbb{R}^n$ which do not belong to $\mathcal{F}$. Meanwhile, observe that for a brick-infinite algebra $A$, there is \emph{a priori} no guarantee that there is a rational $g$-vector outside the $\tau$-tilting fan $\mathcal{F}$. However, this is expected to be true (see Conjecture \ref{Demonet Conj intro}).
\end{remark}

\subsection{Wall-and-chamber structure, canonical decomposition}
To better understand the fan $\mathcal{F}$, we briefly recall some of the ingredients needed to define the wall-and-chamber structure of an algebra, as described in \cite{As2}. 

\medskip
Let us briefly recall that in \cite{As2}, it is shown that for certain bricks $B$ in $\brick(A)$, the set
\[\Theta_B:=\{v \in K_0(\proj(A))_\mathbb{Q} \mid B \text{ is $v$-semistable}\}\]
has codimension one in $\mathbb{Q}^n$ and is called a \emph{wall}. Those walls are naturally seen in the ambient space $\mathbb{R}^n$ (also, see \cite{BST}). The connected components of $\mathcal{F}^\circ$ are called \emph{chambers}. In \cite{As2}, it was proven that any integral $g$-vector is either in a chamber or on a wall. Of course, any two adjacent chambers $C(M,P)$ and $C(M', P')$ are separated by a wall $\Theta_B$. 
 In fact, $B$ is the brick that labels the edge between $\Fac(M)$ and $\Fac(M')$ in the Hasse diagram of the functorially finite torsion classes (see \cite{BCZ} and \cite{DI+}).

For a $g$-vector $v$ in $K_0(\proj(A))$, recall that $\Z(v)$ denotes the $\tau$-regular component corresponding to $v$.
Every such $g$-vector has canonical decomposition, which essentially captures how a map in general position in $\Hom_A(P(v^-), P(v^+))$ decomposes. When a map in general position is injective, that amounts to knowing how a module in general position in $\Z(v)$ decomposes. It happens that these summands are also in general position in other $\tau$-regular components, hence every such summand has a corresponding $g$-vector. In fact, the canonical decomposition of $v$ is a sum of these $g$-vectors (for more details, see \cite{DF} or \cite[Section 9.5]{G+1}). 

For a $g$-vector $v$ in $K_0(\proj(A))$, if $v \in \mathcal{F}^\circ$, then there is a unique basic support $\tau$-tilting pair $(M,P)$ such that $v \in C(M,P)$. That being the case, the canonical decomposition of $v$ is given in terms of the $g$-vectors of the indecomposable summands of $M \oplus P$. For the $g$-vectors outside $\mathcal{F}$, finding the canonical decomposition may be harder. 
Nevertheless, we say $v$ is \emph{indecomposable} if its canonical decomposition is $v$. That means a representation in general position in $\Z(v)$ is indecomposable. 
Some $g$-vectors, even if outside of $\mathcal{F}$, have well-behaved canonical decomposition. A $g$-vector $v$ is said to be \emph{$E$-tame} if the canonical decomposition of $2v$ is $v+v$. 
This is known to be the case if and only if $E(v,v)=0$, where
$$E(v,v):= \min \{{\rm dim}_k\Hom_A(Z,\tau Z')\mid Z, Z' \in \Z(v)\}.$$ 
Note that this minimum is attained on an open set of $\Z(v) \times \Z(v)$, equipped with the Zariski topology
(see also \cite{G+2}).

\section{Bricks and $\tau$-rigid modules}
As before, $\textit{i}\taurigid(A)$ denotes the set of all indecomosable $\tau$-rigid modules in $\modu A$, that is, $\textit{i}\taurigid(A):=\{X \in \ind(A)| \Hom_A(X,\tau_A X)=0\}$. 
We also recall that $\brick(A):= \{M \in \ind(A)| \End_A(M)\simeq k\}$.
In the past few years, primarily motivated by the connections to $\tau$-tilting theory, the interactions between these two sets of indecomposable modules have been studied extensively.
In particular, in \cite{DIJ} the authors established an explicit map, often known as the ``brick-$\tau$-rigid correspondence", from the set of indecomposable $\tau$-rigid modules into the set of bricks.
In the following, we briefly recall the map and some of its main properties used in our work.

\medskip
We begin with a brief recollection of some standard terminology and fix some notations. For details, we only provide references. Throughout, by a subcategory $\mathcal{C}$ of $\modu A$, we always mean $\mathcal{C}$ is a full additive subcategory which is closed under direct sums and direct summands. 
A subcategory $\mathcal{T}$ of $\modu A$ is called a \emph{torsion class} if $\mathcal{T}$ is closed under quotients and extensions.
As shown in \cite{AIR}, the study of $\tau$-tilting theory of $A$ closely relates to that of functorially finite torsion classes in $\modu A$. The latter are the torsion classes of the form $\Fac(M)$, for some $\tau$-rigid modules $M \in \modu A$. Here, $\Fac(M)$ denotes the category of all those modules in $\modu A$ which are quotients of some direct sums of $M$. For details on torsion classes and $\tau$-tilting theory, see \cite{AIR} and \cite{DIJ}.
\medskip

In the following, for $M \in \modu A$,  by $\Lambda_M$ we denote the endomorphism algebra $\End_A(M)$, whose radical is denoted by $\rad_{\Lambda_M}$.

\begin{theorem}[\cite{DIJ}]\label{Thm of DIJ}
Let $A$ be an algebra, and for each indecomposable $\tau$-rigid module $X$, define $B_X:=X/\rad_{\Lambda_X}X$. 
Then, 
\begin{enumerate}
    \item $B_X$ is a brick, and $\Fac(X)$ is the smallest torsion class containing $B_X$.
    
    \item $\Psi: \textit{i}\taurigid(A) \rightarrow \brick(A)$, given by $\Psi(X):=B_X$, is an injective map. 

    \item $\Psi$ is a bijection provided that $\textit{i}\taurigid(A)$ is a finite set.
\end{enumerate}

\end{theorem}

Before stating an important geometric property of the set of bricks specified by the brick-$\tau$-rigid correspondence, we make some helpful observations. In particular, we note if $A$ is hereditary, the map $\Psi$ in the above theorem is the identity map, that is $\Psi(X)=X$. This is because every indecomposable rigid module over a hereditary algebra is known to be a brick. For an arbitrary algebra $A$, if $X$ is an indecomposable $\tau$-rigid module which is not a brick, then the corresponding brick $B_X$ is not $\tau$-rigid (otherwise, $\Psi$ sends $X$ and $B_X$ to the same brick $B_X$, which contradicts its injectivity property). That is to say, $\textit{i}\taurigid(A)\cap \brick(A)$ is the set of elements on which $\Psi$ acts as the identity map, and each $X \in \textit{i}\taurigid(A)\setminus \brick(A)$ is sent to its unique corresponding brick $B_X \in \brick(A) \setminus \textit{i}\taurigid(A)$. 

\medskip
By Theorem \ref{Thm of DIJ}, the image of the map $\Psi$ consists of exactly those $Y \in \brick(A)$ for which the smallest torsion class in $\modu A$ containing $Y$ is functorially finite. These bricks are also used in \cite{BCZ} to give a labelling of the edges of the poset of functorially finite torsion classes in $\modu A$ (also, see \cite{DI+}). 
In our previous work, we showed that such bricks also have an important geometric property. Here we recall this handy property, which be used in the rest of this paper. For the proof and further details, see \cite[Section 6]{MP1}.

\begin{proposition}[\cite{MP1}]\label{Prop:open orbit}
Let $X$ be an indecomposable $\tau$-rigid $A$-module with the corresponding brick $B_X:=X/\rad_{\Lambda_X}X$. Then, $B_X$ has open orbit in the irreducible component that contains it. 
\end{proposition}

As shown in \cite{MP1}, the preceding statement implies that $A$ is brick-infinite if and only if infinitely many $M \in \brick(A)$ have open orbits in their corresponding module variety. On the other hand, the geometric interpretation of the 2nd bBT Conjecture asserts that any brick-infinite algebra should also admit at least one brick (and, therefore, infinitely many non-isomorphic bricks) whose orbit is not open in any module variety containing it (for details, see \cite[Section 6]{Mo1}).

\medskip
The following proposition plays an important role in our work and allows us to associate bricks to arbitrary indecomposable modules.
We particularly thank William Crawley-Boevey for drawing our attention to some earlier results on the study of bricks over hereditary algebras in \cite{Ri2} and \cite{CB1}, and also for generalizing proof of a similar statement in the latter reference. In particular, the following proof is similar to the argument of Lemma 2 in Section 2 of \cite{CB1}. However, for completeness, we provide a full proof and slightly generalize the latter part of the aforementioned argument to handle $\tau$-rigidity. 
Regarding the following result, we particularly note that if $A$ is not hereditary, then in general one cannot replace the condition $\Hom_A(Y,\tau Y)\neq 0$ with the stronger condition $\Ext^1_A(Y,Y) \ne 0$ (similarly, $\Hom_A(\tau^{-} Z, Z)\neq 0$ cannot be replaced with $\Ext^1_A(Z,Z) \ne 0$).

\begin{proposition}\label{Prop:existence of non-tau-rigid bricks}
If $X$ is an indecomposable $A$-module, not a brick, then 
\begin{enumerate}
    \item $X$ admits a proper quotient $Y$ which is a brick and $\Hom_A(Y,\tau Y)\neq 0$.
    \item $X$ admits a proper submodule $Z$ which is a brick and $\Hom_A(\tau^{-} Z, Z)\neq 0$.
\end{enumerate}

\end{proposition}

\begin{proof}
We only prove the first part. The second part follows from a dual argument.

Since $X$ is not a brick, we can choose a nonzero $\phi \in \End_A(X)$ that is not an isomorphism and assume $\Ima(\phi)$ is of minimal dimension. Hence, $\phi$ is nilpotent and the minimality assumption implies $\phi^2=0$.
Consider the decomposition $X/\Ima(\phi)=\bigoplus_{i=1}^t C_i$, where each $C_i$ is a nonzero indecomposable module, and then define $\overline{\phi}:X/\Ima(\phi) \rightarrow \Ima(\phi)$, given by $\overline{\phi}\big(x+\Ima(\phi)\big)=\phi(x)$.
Below, by $\iota:\Ima(\phi)\rightarrow X$ and $\pi:X \to X/\Ima(\phi)$ we respectively denote the natural embedding and projection. Moreover, for $1\leq i \leq t$, we have the natural surjective maps $\pi_i:X/\Ima(\phi)\rightarrow C_i$, and injective maps $\iota_i:C_i \rightarrow X/\Ima(\phi)$. 

Observe that there exists $1 \leq j \leq t$ for which the $\overline{\phi}\iota_j\neq 0$. Fix such a $1 \leq j \leq t$ and define $\alpha:=\overline{\phi}  \iota_j$.
Then, consider the map $\beta:=\iota \alpha \pi_j   \pi$, which clearly belongs to $\End_A(X)$. 
Note that $\beta\neq 0$, and the minimality assumption on $\phi$ implies that $\alpha: C_j \rightarrow \Ima(\phi)$ must be surjective.
We aim to show that $C_j$ is not $\tau$-rigid, which is equivalent to proving that $\Ext^1_A(C_j,\Fac(C_j))\neq 0$. 
We consider the following diagram, where the bottom row is the non-split short exact sequence defined naturally, and the top row is obtained via the pullback of the embedding $\iota_j:C_j \rightarrow \bigoplus_{i=1}^t C_i$. In particular, $\iota_j  b= \pi   c$.
     $$\xymatrix{0 \ar[r] & \Ima(\phi) \ar@{=}[d]  \ar[r]^a& U \ar[d]^{c} \ar[r]^b& C_j \ar[d]^{\iota_j} \ar[r] & 0 \\ 0 \ar[r]& \Ima(\phi) \ar[r]_{\iota} & X \ar[r]_{{\pi}}  & \bigoplus_{i=1}^t C_i \ar[r] \ar[r] & 0}$$

If we assume $\Ext^1(C_j,\Ima(\phi))=0$, the short exact sequence in the top row splits, meaning that $b$ is a split epimorphism and $a$ is a split monomorphism, that is, there exist $s:C_j \rightarrow U$ and $r:U \rightarrow \Ima(\phi)$ such that $b  s=1_{C_j}$ and $r   a=1_{\Ima(\phi)}$.
Consequently, we have $b=(\pi_j \iota_j)  b=(\pi_j \iota_j) (b  s)  b=\pi_j  (\iota_j  b)  s  b=\pi_j  (\pi   c)  s  b$.

This implies $\Big(1_{C_j}- \pi_j\pi c s\Big) b=0$, and because $b$ is a surjective map, we get $1_{C_j}= \pi_j  {\pi}   c  s$. Hence, $c  s$ is a split monomorphism, thus $X=C_j$. This is the desired contradiction. Hence, $\Ext^1(C_j,\Ima(\phi))\neq 0$ and thus $\Ext^1(C_j,\Fac(C_j))\neq 0$.

We have shown that $C_j$ is an indecomposable proper quotient of $X$ which is not $\tau$-rigid.
If $C_j$ is a brick, we are done. Otherwise, we iterate the above argument until we get a proper quotient of $X$, say $Y$, which is a brick and not $\tau$-rigid. 
\end{proof}

We note that, unlike the brick-$\tau$-rigid correspondence in Theorem \ref{Thm of DIJ}, the brick described in the preceding proposition is not necessarily unique. Namely, starting from an arbitrary indecomposable $X$ which is not a brick, one can \emph{a priori} find several non-isomorphic quotients of $X$ which are bricks and not $\tau$-rigid.
Moreover, we note that several non-isomorphic indecomposable $A$-modules that are not bricks may have the same quotient module which is a brick but not $\tau$-rigid. This is, for instance, the case for local algebras, where the only brick is the simple module.

\section{Locally representation-directed algebras}\label{Section:new char. of loc-rep-directed}

As before, let $A=kQ/I$ be an algebra of rank $n$, with the fix $\{e_1, \cdots, e_n\}$ set of primitive idempotents associated to the vertex set $Q_0$. For $M$ and $N$ in $\ind A$, let $\rad_A(M,N)$ denote the set of non-isomorphisms in $\Hom_A(M,N)$.
For some $t \in \mathbb{Z}_{>0}$, and a set of indecomposable modules $M_j \in \modu A$, a sequence $(M_0, ...,M_t)$ is called a \emph{cycle} if $M_0=M_t$ and $\rad_A(M_{j-1},M_j) \neq 0$, for all $1 \leq j \leq t$. An algebra $A$ is said to be a \emph{representation-directed} (rep-directed, for short) if there exists no cycles in $\modu A$.
For a fixed vertex $i \in Q_0$, a cycle $(M_0, ...,M_t)$ in $\modu A$ is said to be an \emph{$i$-cycle} provided $e_i\rad_A(M_{j-1},M_j) \neq 0$, for all $1 \leq j \leq t$. 
Consequently, $A$ is called \emph{locally representation-directed} provided there exists no $i$-cycle in $\modu A$, for all $i \in Q_0$. 
Each directed algebra is evidently locally directed, but the converse is not necessarily true (see part $(1)$ of Example \ref{Examples}).
\medskip

Before we state some properties of the locally rep-directed algebras, let us recall that every $M \in \modu A$ can be viewed as a representation of $(Q,I)$, formally presented by $(\{M_i\}_{i\in Q_0}, \{M(\alpha)\}_{\alpha \in Q_1})$. 
In particular, each $f \in \Hom_A(M,N)$ is of the form $(f_i)_{i\in Q_0}$, where $f_i:M_i\rightarrow N_i$ are some $k$-linear maps between the vector spaces $M_i$ and $N_i$ which satisfy some compatibility conditions. 

We note that several characterizations of locally rep-directed algebras have already appeared in the literature. In particular, from the work Dräxler \cite{Dra}, it is known that an algebra is locally rep-directed if and only if every indecomposable module is a brick (for the proof and other characterizations of these algebras, see \cite{Dra}). 
In the next proposition, we only state some of the results shown by Dräxler, and for the sake of completeness, we provide a proof of this statement.

\begin{proposition}[\cite{Dra}]\label{Prop: local-rep-directe}
If $A$ is locally representation-directed, then every indecomposable $A$-module is a brick and also $\tau$-rigid. In particular, locally representation-directed algebras are representation-finite.
\end{proposition}
\begin{proof}
To show that every $X \in \ind(A)$ is a brick, assume otherwise and take a non-zero map $f \in \rad_A(X,X)$. We have $f=(f_i)_{i\in Q_0}$, and obviously, $f_j:X_j \to X_j$ is non-zero, for some $j \in Q_0$. Hence, $e_j \rad_A(X,X) \neq 0$, which implies that $(X,X)$ is a $j$-cycle, contradicting the assumption that $A$ is locally rep-directed.

To show every indecomposable is $\tau$-rigid, assume otherwise and take $X \in \ind(A)$ with $\Hom_A(X,\tau X)\neq 0$.
In particular, take a non-zero element $h \in \rad_A(X, \tau X)$, where $h=(h_i)_{i\in Q_0}$, hence $h_j\neq 0$, for some $j \in Q_0$.
Because $X$ is not projective, consider the Auslander-Reiten sequence 
$0\rightarrow \tau X \xrightarrow{f} Y \xrightarrow{g} X \rightarrow 0$.
Since $f$ is a monomorphism, $\phi_j:=f_j  h_j: X_j\rightarrow Y_j$ is non-zero. 
Obviously, projection of $\Ima(\phi_j)$ onto one the indecomposable summands of $Y$ is non-zero. Let $Z$ be such an indecomposable summand, with the corresponding projection map $\pi:Y \to Z$. Then, observe that $\pi_j   f_j \neq 0$. Now, we have an irreducible morphism $g':Y \to X$, and because $X$ is supported at vertex $j$, then $g'_j$ is non-zero, due to the irreducibility. This leads to a $j$-cycle $(X, \tau X, Y, X)$. This contradiction the assumption and finishes the proof.

The last assertion follows from the Second Brauer-Thrall theorem and a standard geometric argument. In particular, we note that every ($\tau$-)rigid module has an open orbit in the irreducible component it belongs to.
\end{proof}

\begin{remark}
    We note that if almost all indecomposable modules are bricks (resp. rigid), then $A$ is necessarily representation-finite. 
\end{remark} 
\begin{proof}
     Assume that $A$ is a representation-infinite algebra. Then, by the well-known dichotomy theorem of Drozd, $A$ is tame representation-infinite or wild. If $A$ is tame, then it follows from Crawley-Boevey results (see \cite{CB2}) that the Auslander-Reiten quiver of $\modu A$ contains infinitely many tubes, hence there are infinitely many indecomposable modules that are not brick. If $A$ is wild, then there is a faithful embedding $\modu B \to \modu A$ where $B$ is some tame representation-infinite algebra. Any infinite family of non-isomorphic indecomposable non-bricks over $B$ will remain such a family over $A$. This shows the claim for bricks. 
     
     Now we treat the other case, which is a well known fact. Again, assume that $A$ is representation-infinite. By the 2nd Brauer-Thrall Conjecture (now theorem), there exists an irreducible component $\Z$ in $\Irr(A)$ containing infinitely many orbits of indecomposable modules. By Voigt's Lemma, every rigid module has an open (and dense) orbit, hence there are infinitely many (non-isomorphic) indecomposable modules in $\Z$ which are not rigid.  
\end{proof}
 
The above proposition motivates our new treatment of locally representation-directed algebras from the viewpoint of $\tau$-tilting theory, as well as a geometric generalization of this family considered in Section \ref{Section:geometric generalization of loc-rep-directed}.
In particular, as the main theorem of this section, we show the following result. Below, we use the notations and terminology introduced in the previous sections.

\begin{theorem}\label{Thm:new characterization of locally rep-directed}
For an algebra $A$, the following are equivalent:
\begin{enumerate}
    \item $A$ is locally representation-directed;
    \item Every brick is $\tau$-rigid.
\end{enumerate}
Consequently, if any two of the three sets $\ind(A)$, $\brick(A)$, and $\textit{i}\taurigid(A)$ are equal, then we have $\ind(A)=\brick(A)=\textit{i}\taurigid(A)$.
\end{theorem}

Before proving the above theorem, we make some observations that put this result in perspective.
In particular, regarding part $(2)$, we observe that one cannot replace the notion of $\tau$-rigidity with the rigidity. That is, there exist algebras for which every brick is rigid, but they are not locally rep-directed (see Example \ref{Examples}). We also note that, by the 2nd bBT Conjecture and a standard geometric argument, an algebra for which every brick is rigid is expected to be brick-finite, but we do not know of any proof of this statement in full generality. Some new results in this direction will be shown in the following sections.

\begin{proof}[Proof of Theorem \ref{Thm:new characterization of locally rep-directed}]
The implication $(1) \rightarrow (2)$ follows from Proposition \ref{Prop: local-rep-directe}.
To prove $(2) \rightarrow (1)$, by the main result of \cite{Dra}, it suffices to show $\ind(A)=\brick(A)$. If we assume otherwise, there exists $X \in \ind(A)$ which is not a brick. Hence, by Proposition \ref{Prop:existence of non-tau-rigid bricks}, $X$ admits a proper quotient $Y \in \brick(A)$ which is not $\tau$-rigid. This gives the desired contradiction.

The last assertion is an immediate consequence of the previous arguments. In particular, if every indecomposable is $\tau$-rigid, Proposition \ref{Prop:existence of non-tau-rigid bricks} (or alternatively, the injectivity of the map $\Psi$ in Theorem \ref{Thm of DIJ}) implies $\ind(A)=\brick(A)=\textit{i}\taurigid(A)$.
\end{proof} 

We observe that Theorem \ref{Thm:new characterization of locally rep-directed} implies the following characterization of the locally rep-directed algebras in terms of the brick-$\tau$-rigid correspondence, as well as characterization of torsion classes treated in \cite{DIJ}. As in Theorem \ref{Thm of DIJ}, throughout, by $\Psi: \textit{i}\taurigid(A) \rightarrow \brick(A)$ we denote the injective map which sends $X\in \textit{i}\taurigid(A)$ to its corresponding brick $B_X \in \brick(A)$.

\begin{corollary}\label{Cor: locally rep-directed}
An algebra $A$ is locally representation-directed if and only if the brick-$\tau$-rigid correspondence $\Psi$ is the identity map (and therefore bijective). This is the case if and only if every torsion class in $\modu A$ is of the form $\Fac(X_1\oplus \cdots \oplus X_d)$, for a set of bricks $X_1, \dots, X_d$ satisfying $\Hom(X_i, \tau X_j)=0$, for all $1 \leq i, j \leq d$.
\end{corollary}

Observe that in the first part of the above corollary, the assumptions that $\Psi$ is the identity map is a necessary condition. In particular, one can find algebras for which $\Psi$ is bijective, or algebras for which $\Psi$ acts like the identity map (with a larger codomain), but the algebra is not locally rep-directed (see Example \ref{Examples}).
Before presenting some explicit examples, let us remark that every quotient of a locally rep-directed algebra is again locally rep-directed, and the first non-trivial examples of locally rep-directed algebras should be of rank at least $2$. In fact, if $A=kQ/I$ is a locally rep-directed algebra, where $I$ is an admissible ideal, then $Q$ contains no loops (i.e, no arrow starts and ends at the same vertex), and for any $i$ and $j$ in $Q_0$, there could be at most one arrow from $i$ to $j$.

\begin{example}\label{Examples} 
The following examples clarify some aspects of the problems discussed above, and also motivate the problems considered in the following sections.
\begin{enumerate}
    \item Let $Q$ be the quiver 
\begin{center}
\begin{tikzpicture}
    \node at (-0.1,-0.1) {$\bullet$};
 \draw [->] (0,0) --(1,0);
    \node at (0.5,0.3) {$\alpha$};
    \node at (1.1,-0.1) {$\bullet$};
 \draw [<-] (0,-0.2) --(1,-0.2);
    \node at (0.5,-0.5) {$\beta$};
\end{tikzpicture}
\end{center}
and consider the ideals $I=\langle \alpha\beta, \beta\alpha \rangle$ and $J=\langle \alpha\beta\alpha, \beta\alpha\beta \rangle$ in $kQ$. It is easy to check that $kQ/I$ and $kQ/J$ are rep-finite, but neither of these algebras is rep-directed.
Through an easy computation one can show that $kQ/I$ is locally rep-directed. In contrast, $kQ/J$ is not locally rep-directed, because it admits some bricks which are not $\tau$-rigid. 
Meanwhile, as remarked in \cite{Dra}, we note that every indecomposable $kQ/J$-module is rigid.  More generally, it is not hard to see that for a Nakayama algebra $A$ of rank $n$, locally representation directedness is equivalent to all projectives being of length at most $n$.

\medskip

\item 
Let $Q$ be a Dynkin quiver and by $\Pi(Q)$, denote the corresponding preprojective algebra of $Q$ (for the explicit construction and details on preprojective algebras, see \cite{Ri2} and the references therein). 
Then, $\Pi(Q)$ is a finite dimensional algebra and, except for very small ranks, $\Pi(Q)$ is known to be wild (and in particular rep-infinite).
As shown in \cite{Mi}, the algebra $\Pi(Q)$ is brick-finite.
Moreover, from \cite[Lemma 1]{CB3} it follows that every brick in $\modu \Pi(Q)$ is rigid (also, see \cite{DI+}). However, one can find bricks that are not $\tau$-rigid. 
For instance, if $Q$ is the Dynkin quiver of $A_5$, consider the representation $M$ of $\Pi(Q)$ specified by:
\begin{center}
\begin{tikzpicture}
    \node at (-0.15,-0.1) {$k$};
 \draw [->] (0,0) --(1,0);
    \node at (0.5,0.3) {$0$};
    \node at (1.15,-0.1) {$k$};
 \draw [<-] (0,-0.2) --(1,-0.2);
    \node at (0.5,-0.5) {$1$};
\draw [->] (1.3,0) --(2.3,0);
    \node at (1.75,0.3) {$1$};
    \node at (2.45,-0.1) {$k$};
 \draw [<-] (1.3,-0.2) --(2.3,-0.2);
    \node at (1.75,-0.5) {$0$};
\draw [->] (2.6,0) --(3.6,0);
    \node at (3.75,-0.1) {$0$};
 \draw [<-] (2.6,-0.2) --(3.6,-0.2);
\draw [->] (3.9,0) --(4.9,0);
    \node at (5.05,-0.1) {$0$};
 \draw [<-] (3.9,-0.2) --(4.9,-0.2);
\end{tikzpicture}
\end{center}
where the unlabeled arrows are evidently zero maps.
Then, $M$ is clearly a brick, but a simple computation shows that $M$ is not $\tau$-rigid in $\modu \Pi(Q)$.
Observe that, for $\Pi(Q)$, the brick-$\tau$-rigid correspondence is a bijection, but evidently does not act like the identity map (compare Corollary \ref{Cor: locally rep-directed}). 
\end{enumerate}
\end{example}

To end this section, we mention the following result, which gives an explicit characterization of locally representation-directed gentle algebras. In \cite{GMM}, algebras for which all indecomposable modules are bricks were called brick algebras (so brick algebras are exactly the locally representation-directed algebras). Let $A = kQ/I$ where $I$ is admissible. If $w = \alpha_r^{\epsilon_r}\cdots\alpha_1^{\epsilon_1}$ is a walk in $Q$, where the $\alpha_i$ are arrows and $\epsilon_i =\pm 1$, then we say that $w$ contains a relation from $I$ if there exists $r \in I$ such that $r$ or its formal inverse appears in $w$.

\begin{proposition}[\cite{GMM}] \label{PropGMM}
    Let $A$ be a gentle algebra with $A = kQ/I$ (where $I$ is admissible and generated by quadratic monomial relations). Then $A$ is locally representation-directed if and only if every cyclic walk in $Q$ contains at least two relations from $I$.
\end{proposition}

Note that the characterization of locally representation-directed algebras given in the above proposition does not extend naturally to arbitrary string algebras. As a non-example, let $A$ be a Nakayama algebra given by a cycle-quiver with $n \ge 3$ vertices with the admissible ideal $I$ of relation given by all paths of length $n$. Then, $A$ is a string algebra but is not gentle, and it is easily checked that $A$ is locally representation-directed.

\section{$\tau$-regular components}\label{Section:gen tau-reduced components}
As explained below, the $\tau$-regular components can be treated as a conceptual generalization of $\tau$-rigid modules. While doing so, we translate and generalize some algebraic properties of $\tau$-rigid modules in the geometric setting. Here, we only recall the definition and some properties used in our work, and for more details on $\tau$-regular components, we refer to \cite{BS} and the references therein.

\medskip

As in Section \ref{Subsection: Generically tau-reduced components}, for $\mathcal{Z} \in \Irr(A)$, let $c(\Z):=\min \{\dim(\Z) - \dim(\mathcal{O}_Z) \,| Z \in \Z \}$, which denotes the generic number of parameters of $\Z$ (equivalently, the maximal codimension of an orbit in $\Z$).
Also recall that $h(\mathcal{Z})$ denotes the minimal value of $\Hom_A(Z, \tau Z)$, where $Z \in \Z$. Note that $c(\Z)$ and $h(\mathcal{Z})$ are generically attained (simultaneously) on an open and dense subset of $\Z$, and $c(\Z) \le h(Z)$.
\medskip

Following \cite{BS} and the new terminology proposed therein, we say an irreducible component $\Z \in \Irr(A)$ is \emph{(generically) $\tau$-regular} if $c(\Z)= h(\Z)$. In the literature, these components have also been called \emph{dually strongly reduced} (for instance, in \cite{GLS, Pl}), as well as, \emph{generically $\tau$-reduced} components (for instance, in \cite{G+1, G+2}). 
Recall that $\Irr^{\tau}(A)$ denotes the set of all $\tau$-regular components from $\Irr(A)$.
Observe that those $\Z \in \Irr^{\tau}(A)$ with $c(\Z)=0$ are in bijection with the isomorphism classes of $\tau$-rigid modules, and we simply refer to them as \emph{$\tau$-rigid} components. In general, however, a brick-infinite algebra can admit many $\tau$-regular components that are not $\tau$-rigid. It is expected that every brick-infinite algebra possesses a component $\Z \in \Irr^{\tau}(A)$ with $c(\Z)>0$.
As further discussed in Section \ref{section:constructing rational}, this is a reformulation of the open question posed by Demonet \cite[Question 3.49]{De}.

\medskip

The next important result from \cite{Pl} describes all $\tau$-regular components and gives a surjective map from (integral) $g$-vectors to these components, where they are called ``dually strongly reduced". For a more recent refinement of this result, see \cite{BS}.

\begin{theorem}[\cite{Pl}]\label{Plamondon's Thm}
Let $\Z$ be a $\tau$-regular component of algebra $A$. Then there is a $g$-vector $v$ and an open dense subset $\mathcal{O}$ in $\Hom_A(P(v^-), P(v^+))$ such that the cokernel of all maps in $\mathcal{O}$ belong to $\Z$ and form an open dense subset. Conversely, if $v$ is any $g$-vector, then there is an open dense set $\mathcal{U}$ of $\Hom_A(P(v^-), P(v^+))$ such that the closure of the set $\{{\rm coker}(f) \mid f \in \mathcal{U}\}$ is a $\tau$-regular component.
\end{theorem}

Given a $g$-vector $v$, throughout $\Z(v)$ denotes the corresponding $\tau$-regular component given in Theorem \ref{Plamondon's Thm}. 

\medskip

The following proposition is useful for the reductions in our treatment of some open conjectures (see Section \ref{Section:min-brick-inf}). In particular, it shows how some $\tau$-regular components for algebra $A$ may induce $\tau$-regular components for a quotient algebra $A/J$. Before we state the result, we remark that the algebraic version of this statement is well known. More specifically, if a $\tau$-rigid $A$-module $Z$ is annihilated by a two-sided ideal $J$ in $A$, then $Z$ remains $\tau$-rigid over the quotient algebra $A/J$. 
Hence, the following proposition gives a geometric analogue and generalization of the aforementioned algebraic fact. In particular, it evidently applies to $\tau$-rigid components.

\begin{proposition} \label{GenReduced}
Let $\Z$ be a $\tau$-regular component of an algebra $A$ with associated $g$-vector $v$. If $\Z$ is annihilated by a two-sided ideal $J$ in $A$, then it can be seen as a $\tau$-regular irreducible component for the algebra $A/J$. If $J$ is radical, then the $g$-vector for the quotient can be identified with $v$.
\end{proposition}

\begin{proof}
Under the above assumptions, it is a standard fact that $\Z$ is an irreducible component for $A/J$, so we only show $\Z \in \Irr^{\tau}(A/J)$. By Theorem \ref{Plamondon's Thm}, there is a $g$-vector $v$ in $K_0(K^b({\rm proj}(A)))$ and an open dense subset $\mathcal{O}$ of $\Hom_A(P(v^-), P(v^+))$ such that $\Z$ has an open dense subset $\mathcal{U}$ and each element in $\mathcal{U}$ is a cokernel of an element in $\mathcal{O}$. Consider the space $\Hom_{A/J}(P(v^-)/JP(v^-), P(v^+)/JP(v^+))$. We note that this also corresponds to the $g$-vector $v$, but over the quotient algebra $A/J$ when $J$ is radical. If $J$ is not radical, then $JP(v^-)$ may not lie in the radical of $P(v^-)$ and the $g$-vector for the quotient can be identified with the vector obtained from $v$ by deleting the entries corresponding to the projective modules $Ae_i \subseteq J$. Now, we claim that there is an open set $\mathcal{O}'$ of  $\Hom_{A/J}(P(v^-)/JP(v^-), P(v^+)/JP(v^+))$ such that the cokernels of elements in $\mathcal{O}'$ cover all of $\mathcal{U}$.

In order to show the claim, we let $\pi^+: P(v^+) \to P(v^+)/JP(v^+)$ and also $\pi^-: P(v^-) \to P(v^-)/JP(v^-)$ be the projection maps. 
Furthermore, consider the projection map $\phi$ as $$\Hom_{A}(P(v^-), P(v^+)) \to \Hom_{A/J}(P(v^-)/JP(v^-), P(v^+)/JP(v^+)),$$
given by $f \mapsto \bar f$.

The map $\phi$ is linear and surjective, and therefore it sends open dense sets to open dense sets. More explicitly, $\phi$ is obtained as follows. One considers $\pi^+f$, which vanishes on $JP(v^-)$. Hence, there is a unique map $\bar f$ such that $\pi^+f = \bar f \pi^-$. Set $\mathcal{O}':=\phi(\mathcal{O})$. Let $f \in \mathcal{O}$ and $M:=\coker(f)$. Recall that $M$ is annihilated by $J$. If $g: P(v^+) \to M$ is given by $\coker(f)$, then $P(v^+)J \subseteq {\rm ker}(g)$. That means $g$ factors through $P(v^+)/JP(v^+)$ via $\pi^+$. We can complete to get the following commutative diagram with exact rows:
    $$\xymatrix{P(v^-) \ar[r]^f \ar[d]^{\pi^-} & P(v^+) \ar[r]^g \ar[d]^{\pi^+} & M \ar[r] \ar@{=}[d] & 0 \\ P(v^-)/JP(v^-) \ar[r] & P(v^+)/JP(v^+) \ar[r]  & M \ar[r] & 0}.$$
We observe that the bottom map is $\bar f$. Clearly, the cokernel of $\bar f$ is $M$, and this proves the proposition.\end{proof}

\medskip

For $\mathcal{Z}$ in $\Irr(A)$, recall that we say $\Z$ is \emph{faithful} if there is no non-zero two sided ideal $J$ that annihilates all modules in $\mathcal{Z}$. This does not necessarily imply that $\Z$ contains a faithful module. In the latter case, we say that $\Z$ is a \emph{strongly faithful} component. Evidently, a strongly faithful component is faithful.

\begin{proposition} \label{prop: faithful is open}
For each irreducible component $\Z$ of $A$, the faithful modules in $\Z$ form an open subset.
\end{proposition}

\begin{proof}
Recall that $A=kQ/I$ and let $\Z \in \Irr(A,\textbf{d})$, where $\textbf{d}=(d_i)_{i \in Q_0}$. For each pair $(i,j)$ of vertices of $Q$, consider a set $S_{i,j}$ of paths in $Q$ from $i$ to $j$ that forms a basis of $e_jAe_i$. For a module $M$, and a pair $(i,j)$ of vertices with $e_jAe_i$ non-zero, consider the matrices $M(p)$ for $p \in S_{i,j}$. This matrix is of size $d_j \times d_i$. To it, we associate a row vector $M'(p)$ with $d_id_j$ entries which is just a vector containing the rows of $M(p)$ in the order they appear in the matrix. Now, we create a matrix $M(i,j)$ with $|S_{i,j}|$ rows consisting of all the $M'(p)$ for $p \in S_{i,j}$. We observe that the rank of this matrix is strictly less than $|S_{i,j}|$ precisely when $M$ is annihilated by a non-trivial linear combination of elements in $S_{i,j}$, thus by a non-zero element in $A$. Hence, we can see that the modules annihilated by a non-zero element of $e_jAe_i$ are defined by the modules $M$ with all $|S_{i,j}|$-minors of $M(i,j)$ being zero. This is a closed subset of $\mathcal{Z}$. Since there are finitely many pairs of vertices of $Q$, we see that the non-faithful modules form a closed set of $\mathcal{Z}$.
\end{proof}

While a faithful component may not be strongly faithful, starting from any faithful component $\Z$, we can construct a faithful module induced by $\Z$.

\begin{proposition} \label{chain_faithful}
    Let $\mathcal{Z}$ be a faithful irreducible component. Then there exist some $Z_1, \ldots, Z_r$ in general position in $\mathcal{Z}$ such that $\oplus Z_i$ is faithful.
\end{proposition}

\begin{proof}
We pick $Z_1$ in general position. If $Z_1$ is faithful, we are done. Otherwise, the annihilator ${\rm ann}(Z_1)$ is non-zero. We note that the modules in $\Z$ that are annihilated by ${\rm ann}(Z_1)$ form a closed set in $\Z$. Since $\Z$ is faithful, there exists a module $Z_2 \in \Z$ in general position that does not belong to this closed subset, and therefore ${\rm ann}(Z_1) \not\subseteq {\rm ann}(Z_2)$. This means that ${\rm ann} (Z_1 \oplus Z_2) \subsetneq  {\rm ann}(Z_1)$.
In general, if $Z_1 \oplus Z_2 \oplus \cdots \oplus Z_i$ is not faithful, then there is $Z_{i+1}$ in general position with annihilator ${\rm ann} (Z_{i+1})$ not containing ${\rm ann} (Z_1 \oplus \cdots \oplus Z_i)$ and we get a proper inclusion ${\rm ann} (Z_1 \oplus \cdots \oplus Z_i \oplus Z_{i+1}) \subsetneq {\rm ann} (Z_1 \oplus \cdots \oplus Z_i)$. Eventually, the chain
$$ \cdots \subsetneq {\rm ann} (Z_1 \oplus \cdots \oplus Z_i) \subsetneq \cdots \subsetneq {\rm ann} (Z_1 \oplus Z_2) \subsetneq  {\rm ann}(Z_1)$$
has to stabilize to the zero ideal. 
\end{proof}

\medskip

The next proposition gives a useful knowledge of the projective dimension of the faithful modules that are in general position of irreducible components in $\Irr^{\tau}(A)$.
Let us recall that, for a $g$-vector $v$, by $\Z(v)$ we denote the $\tau$-regular component given by Theorem \ref{Plamondon's Thm}. Furthermore, as in Section \ref{Preliminaries}, by $E(v,v)$ we denote the minimum value of
$\{{\rm dim}_k\Hom_A(Z,\tau Z')\mid Z, Z' \in \Z(v)\}$.

\medskip

It is well known that if $Z$ is a faithful $\tau$-rigid module, then $\pd(Z)\leq 1$ (see \cite[VIII, Lemma 5.1]{ASS}). The following proposition gives a geometric analogue and generalization of this algebraic fact. In particular, the assumptions of the following statement evidently holds for any faithful $\tau$-rigid component.

\begin{proposition} \label{Reduced_pd1}
Let $\Z(v)$ be a faithful $\tau$-regular component and $Z$ in $\Z(v)$ be in general position.
If $\Z(v)$ is strongly faithful, or if $E(v,v)=0$, then $\pd(Z) \leq 1$. In particular, if $c(\Z(v))\leq 1$, then $\pd(Z) \leq 1$.
\end{proposition}

\begin{proof}
We start with the first case, that is, if $\Z = \Z(v)$ is strongly faithful. In particular, $Z$ is faithful.  By the Auslander-Reiten formula, we have $D\overline \Hom_A(Z, \tau Z) \cong \Ext^1_A(Z,Z)$. Moreover, because $\Z$ is $\tau$-regular, dimensions of the two vector  spaces $\Ext^1_A(Z,Z)$ and $\Hom_A(Z, \tau Z)$ are equal. Therefore, any morphism $Z \to \tau Z$ which factors through an injective module must be zero. Since $Z$ is faithful, $DA$ is generated by $Z$, meaning that there is an epimorphism $g: Z^r \to DA$, for some positive integer $r$. If we apply $\Hom(-, \tau Z)$ to $g$, we get a monomorphism 
        $$0 \to \Hom_A(DA, \tau Z) \to \Hom(Z^r, \tau Z).$$
        If $f: DA \to \tau Z$ is a morphism, then $fg: Z^r \to \tau Z$ factors through an injective module, so $fg=0$. This implies that $f=0$, and consequently $\Hom_A(DA, \tau Z)=0$. Hence, the projective dimension of $Z$ is at most one (see \cite[Lemma IV 2.7]{ASS}).

For the second part, it follows from \cite[Theorem 1.5]{G+2} that if $c(\Z)=1$, then $E(v,v)=0$. Thus, we assume that $E(v,v)=0$. 
Hence, for $r \ge 1$, the canonical decomposition of $rv$ is $v + v + \cdots + v$. Therefore, if $Z_1, \ldots, Z_r$ are in general position in $\Z$, then $Z_1 \oplus Z_2 \oplus \cdots \oplus Z_r$ is in general position in the $\tau$-regular component $\Z(rv)$. Using Proposition \ref{chain_faithful}, we see that $\Z(rv)$ has a faithful module $Z_1 \oplus Z_2 \oplus \cdots \oplus Z_r$. It follows from part one that the projective dimension of the latter is at most one, and hence the same holds for its summands $Z_i$.
\end{proof}

Thanks to the preceding proposition, we obtain the following result.

\begin{corollary}\label{Cor:faithful brick comp. on tame}
Let $A$ be a tame algebra and $\Z$ be a brick component of $A$.
If $\Z$ is faithful, then almost all bricks in $\Z$ have projective dimension one and injective dimension one.
\end{corollary}

\begin{proof}
Observe that, because $A$ is tame, we have $c(\Z) \le 1$. Hence, it suffices to consider the case where $c(\Z)=1$. Such a $\Z$ is necessarily $\tau$-regular. This is because, by a result of Crawley-Boevey's \cite{CB2}, a brick $B$ in general position in $\Z$ satisfies $\tau B \cong B$. Hence, $c(\Z)=h(\Z)=1$ (also see \cite{G+1}, and the references therein). 
Thus, provided that $\Z$ is a faithful brick component, by Proposition \ref{Reduced_pd1}, for $Z$ in the general position of $\Z$, we have $\pd(Z)\leq 1$. For the statement on injective dimension, one just considers the opposite (tame) algebra $A^{\rm op}$ and the corresponding (faithful) brick component $\Z^{\rm op}$ obtained using the duality functor.
\end{proof}

The above corollary is particularly interesting in light of the 2nd bBT Conjecture, because the conjecture asserts that any brick-infinite tame algebra must admit a brick component with a $1$-parameter family. Hence, for minimal brick-infinite tame algebras, the conjecture asserts that there is always a faithful such component (for more details, see \cite[Corollary 8.9]{MP4}).
We also note that the faithfulness is a necessary assumption in Corollary \ref{Cor:faithful brick comp. on tame}. For instance, consider the algebra $A=kQ/I$ given by the quiver 
\begin{center}
\begin{tikzpicture}[scale=0.6]
\node at (-10.45,0) {$_1\bullet$};
\draw [->] (-10,0.15) --(-8,0.15);
\node at (-9,0.5) {$\alpha_1$};
\draw [->] (-10,-0.07) --(-8,-0.07);
\node at (-9,-0.5) {$\alpha_2$};
\node at (-7.6,0) {$\bullet_{2}$};
\draw [->] (-7.5,0.15) --(-5.5,0.15);
\node at (-6.5,-0.25) {$\beta$};
\node at (-5.2,0) {$\bullet_{3}$};
\end{tikzpicture}
\end{center}
\noindent
and the ideal $I=\langle \beta \alpha_1, \beta \alpha_2 \rangle$. Observe that $A$ is a string algebra, and thus tame. Also, $A$ admits a $1$-parameter family of bricks in dimension vector $\textbf{d}=(1,1,0)$, induced by the band $\alpha_2^{-1}\alpha_1$. Note that these bricks are of projective dimension $2$, and the brick component containing them is annihilated by the ideal $\langle e_3\rangle$.

\medskip

We end this section with the following remark on the behavior of the $g$-vectors outside the $\tau$-tilting fan and the existence of $1$-parameter families of bricks over tame algebras. This highlights some conceptual connections between several open conjectures (for more details, see \cite[Section 2]{MP3}). In particular, the following result should be known to the experts.

\begin{remark}\label{2nd bBT and Demonet's conj for tame algebras}
For tame algebras, Demonet's Conjecture and 2nd bBT Conjecture are equivalent. 
Namely, if $A$ is tame, the $\tau$-tilting fan is not complete in $K_0(\proj A)_\mathbb{Q}$ if and only if $\brick(A)$ contains an infinite family of bricks of the same dimension.

\begin{proof}
First, observe that for any algebra $A$, the $\tau$-tilting fan $\mathcal{F}_A$ consists of the $g$-vectors $v \in K_0(\proj A)$ for which $c(\Z_v)=0$. Moreover, recall that if $A$ is tame, then $c(\Z) \le 1$, for every indecomposable $\Z \in \Irr(A)$.

Now let $A$ be a tame algebra. From \cite[Theorem 3.8]{PY}, it follows that any indecomposable (integral) $g$-vector outside the $\tau$-tilting fan $\mathcal{F}_A$ yields a $1$-parameter family of bricks (see also \cite{G+1}). Conversely, if $\brick(A)$ contains infinitely many bricks of the same dimension, there exists a brick component $\Z$ with $c(\Z)>0$, thus $c(\Z)=1$. This brick component is necessarily $\tau$-regular, and if $v$ denotes the (integral) $g$-vector corresponding to $\Z$, then $v$ is obviously outside $\mathcal{F}_A$. 
\end{proof}
\end{remark}

In the later parts of this paper, we discuss some other implications between the 2nd bBT Conjecture and Demonet's Conjecture (see Corollary \ref{Cor: E-tame and bricks}).

\section{Minimal brick-infinite algebras}\label{Section:min-brick-inf}

Recall that an algebra $A$ is \emph{minimal brick-infinite} if $\brick(A)$ is infinite, and every proper quotient $A/J$ of $A$ is such that $\brick(A/J)$ is finite. These algebras are studied in \cite{MP1}, mainly from the viewpoint of ($\tau$-)tilting theory, and are called ``minimal $\tau$-tilting infinite algebras". We recall the following property of these algebras, which is used in this paper. As before, we say that a property holds ``almost all" elements of a set $\mathcal{S}$, provided that it holds for all but finite many elements of $\mathcal{S}$. 
 
\begin{proposition}[\cite{MP1}] \label{Prop almost all proj dim one}
Let $A$ be a minimal brick-infinite algebra. Then almost all indecomposable $\tau$-rigid modules are faithful of projective dimension one.
\end{proposition}

\subsection{Cartan matrix as quadratic form}

We fix a Cartan matrix $C_A = (c_{ij})$ for $A$, where $c_{ij}$ is the number of times $S_i$ appears as a composition factor of $P_j = Ae_j$.  Alternatively, we have $c_{ij}= \dim_k(e_iAe_j)$. We let $\langle -,-\rangle$ denote the bilinear form on $K_0(\proj(A))_\mathbb{R}$ given by the Cartan matrix $C_A$.
 Let $T, T'$ be two-term complexes of projectives in degrees $-1,0$. We denote by $g_T, g_{T'}$ the $g$-vectors of $T$ and $T'$, respectively, in $K_0(\proj(A))$. It follows from standard homological considerations that
$$\langle g_T, g_{T'}\rangle =  - {\rm dim}_k\Hom(T,T'[-1]) + {\rm dim}_k\Hom(T,T') - {\rm dim}_k\Hom(T,T'[1]).$$
We let $q$ denote the corresponding quadradic form on $K_0(\proj(A))_\mathbb{R}$ defined by \[q(x) = x^tC_Ax.\]
where $(-)^t$ denotes the transpose of a vector. 
In particular, if $T$ represents a $\tau$-rigid module of projective dimension at most one, then $T$ is isomorphic to $H^0(T)$. That being the case, the above formula  simplifies to $q(g_T) = {\rm dim}_k\Hom(T,T) > 0.$

\begin{lemma} \label{lemma positivity}
    Let $A$ be minimal brick-infinite. Then $q(x)$ is positive on almost all cones of the $\tau$-tilting fan.
\end{lemma}

\begin{proof}
    Since $A$ is minimal brick-infinite, almost all indecomposable $\tau$-rigid modules are faithful of projective dimension one. Hence, if we take a $g$-vector $v$ in the cone corresponding to a faithful (basic) $\tau$-rigid module $T = T_1 \oplus \cdots \oplus T_r$, we get
    \[v = a_1g_1 + \cdots + a_rg_r\]
    where the $a_i$ are positive real numbers and $g_i$ is the $g$-vector of $T_i$. We see that
    \[\begin{array}{ccl}
        q(v) & = & q(a_1 v_1 + \cdots a_r v_r)  \\
    & = & \sum_{i,j}a_ia_j g_i^tC_Ag_j\\
    & = & \sum_{i,j}a_ia_j {\rm dim}_k\Hom(T_i, T_j) - \sum_{i,j}a_ia_j {\rm dim}_k\Hom(T_j, \tau T_i)\\
    & = & \sum_{i,j}a_ia_j {\rm dim}_k\Hom(T_i, T_j)\\
    & > & 0.
    \end{array}\]
\end{proof}

\subsection{The $\tau$-convergence property}\label{Subsection:tau-convergence}

Given a non-zero vector $v = (v_1, \ldots, v_n)$ in $K_0(\proj A)_\mathbb{R}$, we denote by $\bar v$ the vector \[\bar v = \frac{1}{\sum_i|v_i|}v.\]
We observe that if $v = g_T$ is the $g$-vector of a $\tau$-rigid module $T$ with $\pd(T)\leq 1$, then we get
\[q(\bar v) = \frac{1}{\left(\sum_i|v_i|\right)^2}q(v) = \frac{{\rm dim}_k\End(T)}{(\sum_i|v_i|)^2}.\]
Obviously, if $(T_i)_{i \ge 1}$ is an infinite sequence of non-isomorphic $\tau$-rigid modules of projective dimension one, then the denominators in the above expression is unbounded for this family. If the numerators ${\rm dim}_k \End(T_i)$ stay small compared to the denominators, then the sequence $(\overline{g}_{T_i})$ converges to zero. This motivates the following.

\begin{definition}
We say that an algebra $A$ has \emph{the $\tau$-convergence property} if there exists an infinite sequence $(T_j)_{j \ge 1}$ of $\tau$-rigid modules, pairwise non-isomorphic, for which the {\sl infimum} of the rational sequence $\{q(\overline{g}_{T_j})\}_{j \ge 1}$ is zero.
\end{definition}

From the definition, it is obvious that if $A$ has the $\tau$-convergence property, then $A$ is brick-infinite. Moreover, if the intersection of $\brick(A)$ and $\textit{i}\taurigid(A)$ contains infinitely many modules of projective dimension at most one, then $A$ has the $\tau$-convergence property. For instance, for each non-Dynkin quiver $Q$ and the algebra $A=kQ$, we have $\brick(A) \cap \textit{i}\taurigid(A)$ is an infinite set, therefore $A$ has the $\tau$-convergence property. More generally, we have the following remark.

\begin{remark}
    Any representation-infinite algebra with a preprojective component (or a preinjective component) has the $\tau$-convergence property. 
\end{remark}
 
To put the above remark into a better perspective, we note that there also exist brick-infinite algebras with only finitely many $\tau$-rigid bricks, hence no preprojective (nor any preinjective) component. Moreover, note that not every brick-infinite algebra has the $\tau$-convergence property. For instance, take any brick-infinite algebra for which the Cartan matrix $C_A$ is such that $C_A + C_A^t$ is positive definite. Such examples can be easily found in rank $2$. On the other hand, perhaps the $\tau$-convergence property holds for minimal brick-infinite algebras. Consequently, one can ask the following. Observe that an affirmative answer to this question immediately implies that every minimal brick-infinite algebra has the $\tau$-convergence property.

\medskip

 \textbf{Question}:
    Given a minimal brick-infinite algebra and an indecomposable $\tau$-rigid module $T$, can we always bound the dimension of $\End(T)$ by (a constant times) the dimension of $T$?
 \medskip    

Observe that, given an arbitrary indecomposable module $M$, one cannot expect that the dimension of $\End(M)$ be bounded above by a constant times the dimension of $M$. The best that can be done is ${\rm dim}_k \End(M) \le ({\rm dim}_k M)^2$. However, we suspect that better bounds can be obtained for indecomposable $\tau$-rigid modules, at least when $A$ is minimal brick-infinite or tame.
The above property holds for local algebras, and one may hope to prove the general case by induction on the rank. 

\subsection{Constructing rational vectors outside $\mathcal{F}$} \label{section:constructing rational} 

As remarked earlier, it is widely believed that when an algebra is brick-infinite (equivalently $\tau$-tilting infinite), then there has to be a rational vector outside of the $\tau$-tilting fan, or equivalently, the fan $\mathcal{F}_A$ is not complete in $K_0(\proj A)_\mathbb{Q} = \mathbb{Q}^n$. First posed as a question by Demonet, this can be stated in terms of the notion of $E$-finiteness (see Conjecture \ref{Demonet Conj intro}, which we restate below). Recall that an algebra $A$ is called \emph{$E$-finite} if the $\tau$-tilting fan is complete in $K_0(\proj A)_\mathbb{Q}$. Evidently, every brick-finite algebra is $E$-finite. We also note that $E$-finiteness is closed under taking quotients. For more details on the notion of $E$-finiteness and equivalent formulations, we refer to \cite[Section 6]{AI}.

\medskip

\noindent\textbf{Conjecture} (\cite[Question 3.49]{De}):
    Let $A$ be brick-infinite. Then $A$ is not $E$-finite.

\medskip

If there is a rational (or equivalently, integral) vector outside of $\mathcal{F}_A$, from the description of $\tau$-regular components (see Theorem \ref{Plamondon's Thm}) it follows that there is an indecomposable $\tau$-regular component that is not $\tau$-rigid. One can then try to cook up an infinite family of bricks from such a component. In Section \ref{Section:E-tame algebras}, we do that for $E$-tame algebras (see Corollary \ref{Cor: E-tame and bricks}). 
Meanwhile, as explained in the following remark, in the context of Demonet's question, we employ the following reduction argument (compare with Remark \ref{2nd bBT and Demonet's conj for tame algebras}).

\begin{remark}\label{Rem: Reduction of Demonet's Conj}
    To prove the above conjecture, it suffices to verify it for all $g$-tame minimal brick-infinite algebras.
\end{remark}

\begin{proof} 
First, observe that if $A$ is not $g$-tame, then the $\tau$-tilting fan $\mathcal{F}_A$ is not dense in $\mathbb{R}^n$, hence there is a rational $g$-vector outside $\mathcal{F}_A$; equivalently, $A$ is $E$-infinite. 
Hence, assume that $A$ is $g$-tame and consider a minimal brick-infinite quotient $B = A/J$ of $A$, which is also $g$-tame. If the above conjecture holds for $B$ (that is, if $B$ is $E$-infinite), then there is a $g$-vector $\theta$ in $K_0(\proj B)$ which is not rigid (so outside of $\mathcal{F}_B$). It follows from Prop. \ref{GenReduced} that $\theta$ corresponds to a $g$-vector $\theta'$ in $K_0(\proj A)$ and lies outside of $\mathcal{F}_A$. If $J$ is radical, then we can identify $\theta'$ with $\theta$. If $J$ contains a non-zero idempotent, then we identify $\theta'$ as the $g$-vector obtained from $\theta$ by deleting the entries of all $[P_i]$ where the idempotent $e_i$ is in $J$.
\end{proof}

Inspired by the above conjecture and the preceding remark, in this subsection, we provide some tools to produce some rational vectors outside of the $\tau$-tilting fan. In order to do this, we use the $\tau$-convergence property introduced in Section \ref{Subsection:tau-convergence}.

\medskip

Let $(T_i)_{\ge 1}$ be a sequence of non-isomorphic indecomposable $\tau$-rigid modules with $\pd(T_i)=1$. By $v_i$, denote the $g$-vector of $T_i$. We assume that this sequence satisfies the $\tau$-convergence property.
That is, we assume that the {\sl infimum} of $(q(\bar v_i))_{i \ge 1}$ is zero. With no loss of generality, we can select our $T_i$ such that the sequence $(q(\bar v_i))$ converges to zero and the sequence $(\bar v_i)_{i \ge 1}$ converges to a vector $v$.

\begin{lemma} \label{Lemma:localmin}
Let $A$ be $g$-tame  minimal brick-infinite with a sequence $(v_i)_{i\ge 1}$ as above, with $v$ the limit of the $v_i$. Then $v$ is a local minimum of $q$.
\end{lemma}

\begin{proof}
Assume that there is a vector $v'$ near $v$ for which $q$ takes a negative value. 
Let $\{w_i\}_{i \ge 1}$ be a sequence of rational $g$-vectors converging to $v'$. Since we know that $A$ is $g$-tame, we may assume that each rational $g$-vector $w_i$ lies in the $\tau$-tilting fan $\mathcal{F}$. Using Lemma \ref{lemma positivity}, we see that if $w_i$ is close enough to $v'$, and if $v'$ is close enough to $v$, then $q(w_i) > 0$, which is a contradiction to $q(v') < 0$. The value $q(w_i)$ can be negative only for finitely many chambers. When we approach the limit, we will eventually be outside these finitely many chambers, which are all isolated from the limit.
\end{proof}

\begin{lemma} \label{Lemma:localmin_rational}
    For every local minimum $v$ of $q$, there is a sequence $\{u_i\}_{i \ge 1}$ of non-zero rational local minima that converges to $v$.
\end{lemma}

\begin{proof}
    Observe that $q$ is a homogeneous rational quadratic form in $n$ variables, which we denote by $x_1, \ldots, x_n$. Any local minimum $v$ is a zero of all partial derivatives of $q$. Since $q$ is homogeneous and rational, the partial derivatives $\partial q/\partial x_i$ are rational linear forms. We consider the homogeneous system of linear equations given by the vanishing of all $\partial q/\partial x_i$. We denote by $H_q$ its solution space, which is the null-space of a rational matrix and is non-trivial, because $v \in H_q$. We note that $H_q$ has a rational basis. Hence, if $v$ is not itself rational, then we can find the desired sequence.
\end{proof}

The important fact in the proof of Lemma \ref{Lemma:localmin} is that if $v'$ is a $g$-vector that is close enough to $v$ and belongs to the $\tau$-tilting fan $\mathcal{F}$, then $q(v') > 0$. Hence, any $g$-vector $v''$ close enough to $v$ with $q(v'')=0$ must belong to the complement of $\mathcal{F}$ in $\mathbb{R}^n$. Furthermore, Lemma  \ref{Lemma:localmin_rational} guarantees that such a rational $g$-vector $v''$ exists. Hence, we get the following proposition. 

\begin{theorem} \label{thm:rational_g_vector}
    Let $A$ be minimal brick-infinite. If $A$ has the $\tau$-convergence property, then there exists a rational $g$-vector outside of the $\tau$-tilting fan. That is, $A$ is not $E$-finite.
\end{theorem}

\begin{proof}
    Clearly, the result holds if $A$ is not $g$-tame. Hence, we may assume that $A$ is $g$-tame. We let $(T_i)_{i \ge 1}$ be an infinite sequence of non-isomorphic $\tau$-rigid modules, which we can assume all have projective dimension one, since $A$ is minimal brick-infinite. We let $v_i$ denote the $g$-vector of $T_i$. By the $\tau$-convergence property, we may assume that the sequence $(\bar v_i)_{i \ge 1}$ converges and that $\lim q(\bar v_i)=0$. It follows from lemmas \ref{Lemma:localmin} and \ref{Lemma:localmin_rational} that there is a rational vector outside of $\mathcal{F}$.
\end{proof}

We now show the following result. Note that, by Remark \ref{2nd bBT and Demonet's conj for tame algebras}, for a tame algebra $A$, we have that $A$ is not $E$-finite if and only if there exists an infinite family of bricks of the same dimension, that is, the 2nd bBT Conjecture holds for $A$.

\begin{corollary}\label{Cor:tame-almost all bricks tau-rigid}
    Let $A$ be a tame algebra such that almost all bricks are $\tau$-rigid. Then, $A$ is brick-finite.
\end{corollary}

\begin{proof}
    Assume that $A$ is brick-infinite. Let $B$ be a minimal brick-infinite quotient of $A$. Then $B$ is tame and still satisfies the condition that almost all bricks over $B$ are $\tau_B$-rigid. This follows from the well known fact that if $M$ is $\tau_A$-rigid annihilated by an ideal $J$, then it is also $\tau_{A/J}$-rigid. Now, $B$ is tame so it is also $g$-tame by \cite{PY}. Moreover, it has the $\tau$-convergence property since we have an infinite number of $\tau$-rigid bricks. Hence, it follows from Theorem \ref{thm:rational_g_vector} that there is a rational $g$-vector $v$ which is not in $\mathcal{F}$. We can select $v$ to be indecomposable. It follows from \cite{PY} that the irreducible component of $A$ corresponding to the generic cokernels of $P^{v^-} \to P^{v^+}$ forms a $1$-parameter family of bricks. These bricks cannot be $\tau$-rigid, which is a contradiction. This finishes the proof of the corollary.
\end{proof} 

The following proposition highlights the significance of the $\tau$-convergence property in the study of the main conjecture treated in this work. 

\begin{proposition}\label{Prop: 2nd bBT and  tau-convergence for tame algebras}
    The 2ndbBT Conjecture holds for all tame algebras if and only if the $\tau$-convergence property holds for all minimal brick-infinite tame algebras.
\end{proposition}

\begin{proof}
    Assume that the 2ndbBT Conjecture holds for tame algebras. Let $A$ be minimal brick-infinite tame. Then, there exists an irreducible component $\Z$ of $A$ with a $1$-parameter family of bricks. Because $A$ is tame, it follows from Remark \ref{2nd bBT and Demonet's conj for tame algebras} that the component $\mathcal{Z}$ is $\tau$-regular with $c(\Z) > 0$. We let $v$ denote the $g$-vector of $\Z$, which lies outside of the $\tau$-tilting fan $\mathcal{F}_A$. Then observe that $q(v)=0$; this is because for a brick $B$ in general position in $\Z$, we have $B \cong \tau B$, therefore $q(v) = {\rm dim}_k\Hom(B,B) - {\rm dim}_k\Hom(B,\tau B)=0$. Since $A$ is $g$-tame by \cite{PY}, we have a sequence $\{v_i\}_{i \ge 1}$ of rigid $g$-vectors that accumulates to $v$. It follows that $A$ has the $\tau$-convergence property.

    Conversely, assume that every minimal brick-infinite tame algebra has the $\tau$-convergence property. Let $A$ be a tame algebra. To check that the 2nd bBT Conjecture holds for $A$, it suffices to consider a minimal brick-infinite quotient $A'$ of $A$ and prove the 2nd bBT for $A'$. Since $A'$ has the $\tau$-convergence property, it follows from Theorem \ref{thm:rational_g_vector} that there is a rational $g$-vector outside of the $\tau$-tilting fan of $A'$. By Remark \ref{2nd bBT and Demonet's conj for tame algebras}, since $A'$ is tame, this yields a $1$-parameter family of bricks for $A'$.
\end{proof}

\subsection{Faithful $\tau$-regular components} 
In this short subsection, we prove some further consequences of the above results and particularly treat the $\tau$-regular components over minimal brick-infinite algebras. We use the same notations and terminology introduced above.

\begin{proposition} \label{prop: min brick inf. tau-reg comp}
    Let $A$ be a minimal brick-infinite algebra with $\Z$ a $\tau$-regular component. If $\Z$ is not $\tau$-rigid, then $\Z$ is faithful. Moreover, almost all of the $\tau$-rigid components are faithful.
\end{proposition}

\begin{proof}
    Let $\Z$ be $\tau$-regular and not $\tau$-rigid. Let $v$ be the $g$-vector associated to it. If $J$ is a two-sided ideal which annihilates $\Z$, then it follows from Proposition \ref{GenReduced} that $A/J$ has a $\tau$-regular component that is not $\tau$-rigid. This implies that $A/J$ is brick-infinite, which is only possible when $J=0$ by minimality of $A$. The second assertion follows from Proposition \ref{Prop almost all proj dim one}.
\end{proof}

Although a faithful $\tau$-regular component may not give rise to modules of projective dimension one, we have the following result on a family of algebras which are of interest.
We recall that, to prove the 2nd bBT Conjecture in full generality, it suffices to treat those minimal brick-infinite algebras over which almost all bricks are faithful (see \cite{MP1}). Therefore, any further knowledge of such algebras can provide new impetus to the study of this open conjecture. As such, we show the following result.

If $\theta \in K_0(\proj A)_\mathbb{Q}$, then the ray of $\theta$, denoted by $[\theta]$, is defined as
$$[\theta] = \{r\theta \mid r \in \mathbb{Q},\;r\theta \in K_0(\proj A)\}.$$
We say $[\theta]$ is \emph{strongly faithful} if there is $\theta' \in [\theta]$ such that $\Z(\theta')$ is strongly faithful, that is, if $\Z(\theta')$ contains a faithful module. Observe that if $\theta$ is rigid (that is, is the $g$-vector of a $2$-term silting complex), then for $\theta', \theta'' \in [\theta]$, we have $\Z(\theta')$ is (strongly) faithful if and only if $\Z(\theta'')$ is (strongly) faithful.

\begin{proposition}\label{Prop:measure zero}
Let $A$ be minimal brick-infinite with almost all bricks faithful. Then the non-strongly faithful rays forms a subset which is included in a finite union of some chambers and walls.
\end{proposition}

\begin{proof}
Let $S$ denote the set of all walls $\Theta_B$ corresponding to non-faithful bricks. Observe that $S$ is a finite set. If $\theta$ is any $g$-vector that lies neither in a chamber nor on a wall $\Theta_B$ from $S$, then $\theta$ lies on some $\Theta_U$ where $U$ is a faithful brick. Now, a module $Z$ in general position in $\Z(\theta)$ is such that the Schofield semi-invariant $c^Z(U) \ne 0$. It follows from \cite[Theorem 1]{DW} that the projective dimension of $Z$ is at most one. Note that the proof of \cite{DW} assumes that the quiver is acyclic, but the proof works verbatim if the quiver has oriented cycles. Also, since $A$ is minimal brick-infinite, it follows from Proposition \ref{prop: min brick inf. tau-reg comp} that $\Z$ is faithful. Moreover, it follows from Proposition \ref{chain_faithful} that there are modules $Z_1, \ldots, Z_r$ in $\Z$ which are in general position and such that $Z_1 \oplus \cdots \oplus Z_r$ is faithful. Using that $Z_1 \oplus \cdots \oplus Z_r$ has projective dimension at most one, it follows  from \cite[Theorem 6.1]{BS} that $Z_1 \oplus \cdots \oplus Z_r$ belongs to $\Z(r\theta)$, which is strongly faithful. Hence, $[\theta]$ is strongly faithful. 

If $\theta$ lies in a chamber, then a module in general position in $\Z(\theta)$ is $\tau$-rigid. Since $A$ is minimal brick-infinite, almost all $\tau$-rigid modules are faithful. The proposition now follows.
\end{proof}

\section{Geometry of bricks and $\tau$-rigidity }\label{Section:geometric generalization of loc-rep-directed}

We follow our notations from Section \ref{Preliminaries}. In particular, $A=kQ/I$ is an algebra of rank $n$, and $\Irr(A)$ denotes the set of all irreducible components of all representation varieties of $A$. Moreover, recall that $\Irr^{\tau}(A)$ denotes the set of all $\tau$-regular irreducible components, that is, $\Irr^{\tau}(A):= \{\Z \in\Irr(A)|\, c(\Z)=h(\Z)\}$. 
\medskip

We now consider three special subsets of irreducible components that are important in the following. 
In particular, let $\Irr^i(A)$ denote the set of indecomposable components in $\Irr(A)$, and by $\Irr^b(A)$ we denote the set of brick components in $\Irr(A)$. Finally, by $\Irr^{i\tau}(A)$ we denote the set of all indecomposable $\tau$-regular components, namely, $\Irr^{i\tau}(A):= \Irr^{i}(A) \cap \Irr^{\tau}(A)$.
Observe that the indecomposable $\tau$-rigid modules in $\modu A$ are in bijection with the components $\Z$ in $\Irr^{i\tau}(A)$ for which $c(\Z)=0$.
As remarked before, the set $\Irr^{i\tau}(A)$ can be seen as a geometric generalization of the set $\textit{i}\taurigid(A)$.

\medskip

At the level of modules, we have $\brick(A) \subseteq \ind(A)$ and $\textit{i}\taurigid(A)\subseteq \ind(A)$. For the irreducible components, we get $\Irr^b(A)\subseteq \Irr^i(A)$ and $\Irr^{i\tau}(A) \subseteq \Irr^{i}(A)$, where the first containment holds because each brick component is indecomposable, and the second containment is by definition. Inspired by the above analogy and our results in Section \ref{Section:new char. of loc-rep-directed}, the following question arises naturally. 
To put this question into perspective, recall that, among the three sets $\ind(A)$, $\brick(A)$, and $\textit{i}\taurigid(A)$, two of them are equal if and only if all three are the same (see Theorem \ref{Thm:new characterization of locally rep-directed}).
Hence, one can naturally ask if a similar phenomenon holds among the geometric counterparts.

\medskip

\noindent
\textbf{Question:}
Among the three sets $\Irr^i(A)$, $\Irr^b(A)$, and $\Irr^{i\tau}(A)$, if two of them are equal, does it  necessarily imply that all three sets must be the same?
\medskip

For hereditary algebras, it is known that each representation variety $\Z = \rep(Q,\textbf{d})$ is irreducible, and such a $\Z$ is indecomposable if and only if $\textbf{d}$ is a Schur root (that is, $\Z$ is a brick component). Moreover, every irreducible component is $\tau$-regular. Hence, if $A$ is a hereditary algebra, then $\Irr^i(A)=\Irr^b(A)=\Irr^{i\tau}(A)$. However, the converse is not necessarily true (see Examples \ref{ExGluing2}). 

\medskip

We now present the main result of this section. 
On the one hand, motivated by the above question, it treats a family of algebras with some nice conditions on their irreducible components. On the other hand, it leads to some new results on the 2nd bBT Conjecture. 
More specifically, part $(1)$ implies that the family of algebras considered in this theorem can be considered as a natural generalization of the locally representation-directed algebras (see Theorem \ref{Thm:new characterization of locally rep-directed}). In particular, this family simultaneously contains all hereditary algebras and locally representation-directed algebras. Moreover, the gluing construction given in the rest of this section generates many additional examples that belong to this family. 
Furthermore, in light of the natural reduction of the 2nd bBT Conjecture to minimal brick-infinite algebras, part $(2)$ verifies this open conjecture for some new families.

\begin{theorem}\label{Thm: geometric 2 out of 3}
Let $A$ be an algebra for which two of the three sets $\Irr^i(A)$, $\Irr^b(A)$, and $\Irr^{i\tau}(A)$ are equal. Then, $\textit{i}\taurigid(A) \subseteq \brick(A)$, and furthermore we have:
\begin{enumerate}
    \item $A$ is brick-finite if and only if $A$ is locally representation-directed.
    \item If $A$ is  minimal brick-infinite then $A$ admits (at least) one infinite family of bricks in the same dimension.
\end{enumerate}
\end{theorem}

\begin{proof}
If $\Irr^i(A) = \Irr^b(A)$, then using the components in $\Irr^i(A)$ which are orbit closures of indecomposable $\tau$-rigid modules, we get that every indecomposable $\tau$-rigid module is a brick. The same conclusion holds if $\Irr^b(A) = \Irr^{i\tau}(A)$. Then, assume $\Irr^i(A) = \Irr^{i\tau}(A)$. Note that every brick with an open orbit gives rise to a component $\mathcal{Z}$ in $\Irr^i(A)$ with $c(\mathcal{Z})=0$. Hence, by $\Irr^i(A) = \Irr^{i\tau}(A)$, every such component $\mathcal{Z}$ must be $\tau$-regular, and $c(\mathcal{Z})=0$ forces the given brick to be $\tau$-rigid. Recall that in \cite{MP1}, we proved that every left-finite brick has open orbit. Therefore, using the DIJ map recalled in Theorem \ref{Thm of DIJ}, which is always injective, we get that every indecomposable $\tau$-rigid module is a brick. This proves the first assertion of the theorem; namely, $\textit{i}\taurigid(A) \subseteq \brick(A)$.

To show $(1)$, assume $A$ is brick-finite. Thus, the DIJ map is a bijection. Since every indecomposable $\tau$-rigid module is a brick, we get $\textit{i}\taurigid(A) = \brick(A)$, hence $A$ is locally representation directed (see Theorem \ref{Thm:new characterization of locally rep-directed}). The converse is clear.

To show $(2)$, let $A$ be minimal brick-infinite. Since every indecomposable $\tau$-rigid module is a brick, and because almost all indecomposable $\tau$-rigid modules have projective dimension one, the algebra $A$ has the $\tau$-convergence property. Then, Theorem \ref{thm:rational_g_vector} implies that there exists an (indecomposable) $\tau$-regular component $\mathcal{Z}$ with $c(\mathcal{Z}) > 0$. If $\mathcal{Z}$ is a brick component, then we are done. Otherwise, we are in the situation where $\Irr^i(A) = \Irr^{i\tau}(A)$. 
If there is no infinite family of bricks in the same dimension, then every brick has an open orbit. That being the case, from $\Irr^i(A) = \Irr^{i\tau}(A)$, we conclude that every brick is $\tau$-rigid. This implies that $A$ is locally representation directed (see Theorem \ref{Thm:new characterization of locally rep-directed}), and therefore $A$ must be representation-finite. This is the desired contradiction.
\end{proof}

In the rest of this section, we consider an explicit construction that produces concrete non-trivial examples of those algebras considered in Theorem \ref{Thm: geometric 2 out of 3}. 
For this purpose, we recall some handy techniques, based on the notion of gluing, which have been already considered in the literature (see \cite{Ma1, Ma2}, \cite{Ri3}, and \cite{Bo}).
For the most part, we only provide reference to the known results on the algebraic aspect of this construction.
However, we need to make some observations on the geometric aspects.
As in the above-mentioned references, the arguments presented below can be given in terms of the more rigorous language of functors. However, the more elementary approach used below suffices for our purposes here, which is mainly aimed at producing new examples. For the sake of brevity, some straightforward arguments are left to the reader.

\subsection*{Gluing of algebras}
Let $A=kQ/I$ and $B=kQ'/I'$ be a pair of algebras. As before, $I$ and $I'$ are assumed to be admissible ideals, and we let $R$ and $R'$ denote two fixed sets of admissible relations, respectively in $kQ$ and $kQ'$, that generate $I$ and $I'$.
Henceforth, assume $Q_0=\{x_1, \cdots, x_n\}$ and $Q'_0=\{y_1, \cdots, y_m\}$, which respectively denote the vertex set of $A$ and $B$.
Furthermore, assume $``x_n"$ is a sink in $Q$, and $``y_1"$ is a source in $Q'$. Then, by $\Lambda:=k\Tilde{Q}/\Tilde{I}$ we denote the algebra obtained from the \emph{gluing} of $A$ and $B$ at the vertices $``x_n"$ and $``y_1"$, with the bound quiver $(\Tilde{Q}, \Tilde{I})$ constructed as follows, and illustrated in Example \ref{Exmaple of gluing}: 
\medskip

$\bold{\Tilde{Q}}$: The new quiver is obtained from $Q$ and $Q'$, by identifying the vertices $``x_n"$ and $``y_1"$, and keeping the rest of $Q$ and $Q'$ untouched. Thus, the vertex set of $\Tilde{Q}$ is $\Tilde{Q}_0=(Q_0-\{z\})\cupdot(Q'_0-\{a\}) \cupdot \{v\}$, where $``v"$ denotes the vertex obtained from the identification of $``x_n"$ and $``y_1"$. Hence, $\Tilde{Q}_0=\{x_1,\dots,x_{n-1},v,y_2,\dots,d_m\}$, and the arrow set of $\Tilde{Q}$ is $\Tilde{Q}_1= Q_1\cupdot Q'_1$.

\medskip
$\bold{\Tilde{I}}$: The ideal in $k\Tilde{Q}$ is generated by the set $\Tilde{R}:= R \cupdot R' \cupdot R_v$, where $R_v$ denotes the set of all paths of length $2$ in $\Tilde{Q}$ passing through $``v"$. That is, each element of $R_v$ is of the form $\beta\alpha$, for an arrow $\alpha \in Q_1$ ending at $``x_n"$, and an arrow $\beta \in Q'_1$ starting at $``y_1"$.

\medskip
Observe that $(\Tilde{Q}, \Tilde{I})$ contains a copy of $(Q,I)$, as well as a copy of $(Q',I')$, as bound subquivers.
Moreover, beside the vertex $``v"$ in $\Tilde{Q}$, at which $Q$ and $Q'$ are identified, any other vertex in $\Tilde{Q}_0$ belongs either to $Q_0$ or $Q'_0$, and not both. 
We also note that $A$ and $B$ may have several sinks and sources in their bound quivers, and therefore gluing along each pair gives rise to different algebras. Also, it is trivial that $\Lambda=k\Tilde{Q}/\Tilde{I}$ is isomorphic to $A=kQ/I$ if and only if the bound quiver of $B$ consists of only one vertex, that is, if $B\simeq k$.
Before we discuss some general properties of the gluing, let us present an explicit example to clarify the construction.

\begin{example}\label{Exmaple of gluing}
Let $A=kQ/I$ and $B=kQ'/I'$, where $Q$ and $Q'$ are given by
\begin{center}
\begin{tikzpicture}[scale=0.6]
\node at (-11.7,0) {$Q$:};
\node at (-10.45,0) {${x_1}\bullet$};
\draw [->] (-10,0.15) --(-8,0.15);
\node at (-9,0.5) {$\alpha_1$};
\draw [->] (-10,-0.07) --(-8,-0.07);
\node at (-9,-0.5) {$\alpha_2$};
\node at (-7.6,0) {$\bullet {x_2}$};

\node at (-2,0) {$Q'$:};
\node at (-0.5,-2) {$y_1\bullet$};
\draw [->] (-0.25,-1.8) to (-0.25,0);
\node at (-0.55,-1) {$\beta_1$};
\node at (-0.5,0.1) {$y_2\bullet$};
\draw [->] (0,-2)--(1.8,-2);
\node at (1,-2.35) {$\beta_2$};
\node at (2.2,-2) {$\bullet y_3$};
\draw [->] (2,-1.8) --(2,0) ;
\node at (2.25,0.1) {$\bullet y_4$};
\draw [->] (0,0.2) --(1.8,0.2);
\node at (1,0.5) {$\beta_3$};
\node at (2.35,-1) {$\beta_4$};
\end{tikzpicture}
\end{center}
\noindent
and the ideals $I$ and $I'$ are respectively generated by $R=\emptyset$ and $R'=\{\beta_3\beta_1, \beta_4\beta_2\}$.
Observe that $A=kQ/I$ is the Kronecker algebra, and hence it is hereditary, and $B=kQ'/I'$ is a gentle algebra with the two quadratic relations specified in $R'$.

Glue the algebras $A=kQ/I$ and $B=kQ'/I'$ by identifying the sink $``x_2"$ in $Q$ with the source $``y_1"$ in $Q'$. As the result, we obtain the algebra $\Lambda=k\Tilde{Q}/\Tilde{I}$, where  $\Tilde{Q}$ is the following quiver, in which the new vertex $``v"$ is denoted by $\circ$:
\begin{center}
\begin{tikzpicture}[scale=0.6]
\node at (-2.9,-2.1) {${x_1}\bullet$};
\draw [->] (-2.4,-2.1) --(-0.5,-2.1);
\draw [->] (-2.4,-1.9) --(-0.5,-1.9);
\node at (-1.4,-1.6) {$\alpha_1$};
\node at (-1.4,-2.4) {$\alpha_2$};
\node at (-0.2,-2) {$\circ$};
\node at (-0.2,-2.3) {$v$};

\draw [->] (-0.25,-1.8) to (-0.25,0);
\node at (-0.55,-1) {$\beta_1$};
\node at (-0.5,0.1) {$y_2\bullet$};
\draw [->] (0,-2)--(1.8,-2);
\node at (1,-2.35) {$\beta_2$};
\node at (2.2,-2) {$\bullet y_3$};
\draw [->] (2,-1.8) --(2,0) ;
\node at (2.25,0.1) {$\bullet y_4$};
\draw [->] (0,0.2) --(1.8,0.2);
\node at (1,0.5) {$\beta_3$};
\node at (2.35,-1) {$\beta_4$};
\end{tikzpicture}
\end{center}

Here, the admissible ideal $\Tilde{I}$ is generated by the set $\Tilde{R}:= R \cupdot R' \cupdot R_v$, where $R_v$ consists of all paths of length two in $\Tilde{Q}$ passing through vertex $``v"$, that is $R_v=\{\beta_1\alpha_1, \beta_2\alpha_1, \beta_1\alpha_2, \beta_2\alpha_2 \}$. We observe that $A$ is a $4$-dimensional hereditary algebra of rank $2$, and $B$ is an $8$-dimensional gentle algebra of rank $4$, whose global dimension is $2$. The algebra $\Lambda=k\Tilde{Q}/\Tilde{I}$ obtained from the gluing of $A$ and $B$ is an $11$-dimensional string algebra of rank $5$, whose global dimension is $3$.

Finally, we observe that the algebras $A$ and $B$ can be glued through a different choice of vertices, that is, at the sink $y_4$ in $B$ and the source $x_1$ in $A$.
\end{example}

The notations and terminology introduced in the above construction will be preserved throughout the rest of this section.
Note that $A=\Lambda/J_A$ and $B=\Lambda/J_B$, where $J_A$ and $J_B$ are two ideals in $\Lambda$ that can be generated by idempotents.
In particular, if a module $M \in \modu \Lambda$ is supported only on $Q$ (respectively, on $Q'$), then we can naturally regard $M$ as a module in $\modu A$ (respectively, in $\modu B$).

\medskip
The next lemma captures some basic properties of the algebras obtained from the gluing construction. We use the same notation as above and encourage the reader to consult Example \ref{Exmaple of gluing} while reading the lemma. The proof is left to the reader, as it can be concluded from the construction above. Alternatively, one can show these statements using the detailed study of gluing construction in \cite{Ma1, Ma2}.
\begin{lemma}\label{Lem:on gluing}
Let $A=kQ/I$ and $B=kQ'/I'$ be a pair of algebras and suppose $\Lambda=k\Tilde{Q}/\Tilde{I}$ is obtained from a gluing $A$ and $B$. Then, $\Lambda$ is finite dimensional and it is of rank $|A|+|B|-1$. Furthermore, the following hold:

\begin{enumerate}
\item $\Lambda$ is representation-finite if and only if both $A$ and $B$ are representation-finite.

\item If $M$ is indecomposable, then $M$ is either entirely supported over $Q$ or entirely supported over $Q'$.

\item $\Lambda$ is brick-finite if and only if both $A$ and $B$ are brick-finite.

\end{enumerate}
\end{lemma}

Next, we discuss the interactions between the indecomposable components in $\Irr(A)$ and $\Irr(B)$ with those belonging to $\Irr(\Lambda)$. For the simple module $S_v$ corresponding to the node $v \in \Tilde{Q}_0$, we obviously have $\mathcal{O}_{S_v} = \overline{\mathcal{O}}_{S_v}$ is an indecomposable component.
If $\Z$ is an indecomposable irreducible component of $\Lambda$ different from $\mathcal{O}_{S_v}$, by Lemma \ref{Lem:on gluing}, any indecomposable module in $\Z$ is supported entirely over $Q$ or entirely over $Q'$. 
By irreducibility, only one situation occurs. Conversely, each indecomposable irreducible component for $A$ can be seen naturally as an indecomposable irreducible component for $\Lambda$ where the dimension vector is supported entirely on $Q$. A similar statement holds for $B$, and we get the following proposition.

\begin{proposition} We have $\Irr^*(\Lambda) = \mathcal{O}_{S_v} \cup {\Irr^*(A)} \cup \Irr^*(B)$ where $* \in \{i, b, i\tau\}$.
\end{proposition}

\begin{proof}
    The first and second statements follow from our previous analysis. For indecomposable $\tau$-regular components, we need to make use of Proposition \ref{GenReduced}. We also observe that a simple module over a finite dimensional basic $k$-algebra $C$ is $\tau$-rigid if and only if there is no loop in the quiver of $C$ at the vertex corresponding to that simple module. Hence $S_v$ is $\tau$-rigid seen as an $A,B$ or $\Lambda$-module.
\end{proof}

Earlier in this section, we observed that all hereditary algebras, as well as every locally representation-directed algebra, fall into the setting of Theorem \ref{Thm: geometric 2 out of 3}.  Moreover, in Proposition \ref{PropGMM} we recalled the concrete characterization of those gentle algebras which are locally representation-directed, hence satisfying the desired condition of Theorem \ref{Thm: geometric 2 out of 3}.
Thanks to the above observations, we obtain a large pool of new examples via gluing. This, in particular, produces rep-infinite algebras which are neither hereditary nor gentle but still satisfy the assumptions of Theorem \ref{Thm: geometric 2 out of 3}. For all such algebras, the 2nd bBT Conjecture holds. We summarize this construction in the following example.

\begin{example} \label{ExGluing2}
Let $\Lambda$ be obtained from a gluing of the algebras $A$ and $B$.
If $\Irr^*(A)=\Irr^{\Delta}(A)$ and $\Irr^*(B)=\Irr^{\Delta}(B)$, then $\Irr^*(\Lambda)=\Irr^{\Delta}(\Lambda)$, where $*$ and $\Delta$ are any choice of pairs from $\{i, b, i\tau\}$.
In other words, if a ``two out of three" property simultaneously holds for $A$ and $B$, then the same ``two out of three" property also holds for the algebra obtained from their gluing. 
In particular, for each $n \in \mathbb{Z}_{>0}$, there exist infinitely many (non-isomorphic) algebras of rank $n>1$ that are non-hereditary and each one of them, say $\Lambda$, satisfies $\Irr^i(\Lambda)=\Irr^{b}(\Lambda)=\Irr^{i\tau}(\Lambda)$. Furthermore, for the fixed rank $n$, one can construct $\Lambda$ to be rep-finite, or tame rep-infinite, or of wild type. From Theorem \ref{Thm: geometric 2 out of 3}, it follows that if $\Lambda$ is rep-infinite, then it is in fact brick-infinite and admits at least one parameter family of bricks of the same dimension.
\end{example}

\section{Torsion classes of modules varieties and $E$-tameness}\label{Section:E-tame algebras}

In this section, we first consider the set of torsion classes that naturally arise from the irreducible components of the representation varieties of algebras under consideration. Then, we show that brick components and $\tau$-regular components manifest some interesting properties with respect to the finitely generated torsion classes they determine. Finally, with regard to the open conjectures treated in the previous sections, we extend some of our results on tame algebras to the larger family of $E$-tame algebras.

\medskip

\subsection{Torsion classes of irreducible components} 
For an algebra $A$, recall that a subcategory $\mathcal{T}$ of $\modu A$ is said to be a \emph{torsion class} if $\mathcal{T}$ is closed under quotients and extensions. Let $\tors(A)$ denote the set of all torsion classes in $\modu A$. Moreover, for any $X$ in $\modu A$, by $\mathcal{T}_X$ we denote the smallest torsion class in $\modu A$ which contains $X$. More explicitly, $\mathcal{T}_X=\Filt(\Fac(X))$; hence, each $Y$ in $\mathcal{T}_X$ has a filtration $0=Y_0\subset Y_1 \subset \cdots \subset Y_{t-1}\subset Y_t=Y$ such that, for every $1 \leq i \leq t$, there exists an epimorphism $\psi_i:X^{d_i}\twoheadrightarrow Y_{i}/Y_{i-1}$ for some $d_i \in \mathbb{Z}_{>0}$.

\medskip

Following \cite{Ri1}, we say a torsion class $\mathcal{T}\in \tors(A)$ is \emph{finitely generated} if $\mathcal{T}=\mathcal{T}_X$, for some $X\in \modu A$ (obviously, $X$ is not necessarily unique). 
Every functorially finite torsion class is evidently finitely generated, but the converse is not true in general. We also remark that the finitely generated torsion classes are also called \emph{compact} torsion classes. This is because, in the lattice $\tors(A)$, the finitely generated torsion classes are precisely the elements of $\tors(A)$ that are compact (in the lattice-theoretical sense); for more details, see \cite[Prop. 3.2]{DI+}.

\medskip

In \cite[Theorem 2.2]{Ri1}, the author has recently shown the following important characterization of finitely generated torsion classes. Before we state this result, recall that a (possibly infinite) set of pairwise Hom-orthogonal bricks in $\modu A$ is called a \emph{semibrick}. For some recent studies on semibricks and their connections to the problems treated in this work, we refer to \cite{As1}, \cite{MP3}, and references therein.

\begin{thm}[\cite{Ri1}]

\label{Ringel SB Theorem}
For any algebra $A$, there is a bijection $\Psi$ between finitely generated torsion classes and isoclasses of finite semibricks. If $M \in \modu A$, then there is a unique (up to isomorphism) finite semibrick $\Psi(M)$ that is a quotient of $M$ with $\mathcal{T}_M = \mathcal{T}_{\Psi(M)}$.
\end{thm}

For the explicit description of the map $\Psi$ in the preceding theorem, see \cite{Ri1}.
\medskip

In the following, we are particularly interested in the torsion classes induced by each irreducible component. More specifically, for every $\Z \in \Irr(A)$, we set 
$$\tors(\Z):=\{\mathcal{T}_X \mid X \in \Z\},$$
which consists of all finitely generated torsion classes induced by the modules in $\Z$. 

\medskip

The following statement concerns the size of $\tors(\Z)$. This particularly indicates how the study of $\tors(\Z)$ can provide a new perspective to our treatment of some open conjectures considered in the previous sections (see Remarks \ref{Rem: on tors(Z)} and \ref{Rem: stability}).

\begin{proposition}\label{Prop: infinite tors(Z)}
The algebra $A$ admits an infinite family of (non-isomorphic) bricks of the same dimension if and only if $\tors(\Z)$ is an infinite set, for some $\Z \in \Irr(A)$. 
\end{proposition}

\begin{proof}
    The forward direction follows from the fact that for two non-isomorphic bricks $B, B'$, we have $\mathcal{T}_B \ne \mathcal{T}_{B'}$. For the other direction, let $\{Z_i\}_{i\in I}$ be an infinite family of modules in $\Z$ such that the corresponding torsion classes $\{\mathcal{T}_{Z_i}\}_{i\in I}$ are all distinct. Let $d$ denote the dimension of modules in $\Z$. It follows from Theorem \ref{Ringel SB Theorem} that the semi-bricks $\{\Psi(Z_i)\}_{i\in I}$ are pairwise non-isomorphic. Since an indecomposable brick summand of any $\Psi(Z_i)$ has dimension bounded by $d$, this yields infinitely many non-isomorphic bricks of dimension bounded by $d$. Since there are only finitely many irreducible components for modules of dimension at most $d$, the desired conclusion follows.
\end{proof}

\begin{remark}\label{Rem: on tors(Z)}
For a brick-infinite algebra $A$, observe that Proposition \ref{Prop: infinite tors(Z)} implies that the 2nd bBT Conjecture holds if and only if there exists $\Z \in \Irr(A)$ such that $\tors(\Z)$ is an infinite set. As noted in Remark \ref{2nd bBT and Demonet's conj for tame algebras}, this has some immediate applications to some other open conjectures (for details, see \cite[Section 2]{MP3}). 
In our future work, we study $\tors(\Z)$ more extensively and show how $\tors(\Z)$ captures and reflects several interesting geometric properties of $\Z$ in terms of the induced torsion classes.
\end{remark}

In light of the preceding remark, it is interesting to find $\Z \in \Irr(A)$ for which $\tors(\Z)$ is an infinite set. 
Below, we list some specific cases where this property holds. 
We note that a statement equivalent to part (2) of the next proposition is shown independently in \cite{Pf2}. Nevertheless, we give an alternative proof of this result (see also Remark \ref{Rem: stability}).

\begin{proposition} \label{Prop on c(Z) zero}
For any algebra $A$, we have
\begin{enumerate}
    \item If $\Z \in \Irr^b(A)$ with $c(\Z)>0$, then $\tors(\Z)$ is an infinite set. 

    \item If $\Z \in \Irr^{\tau}(A)$ with $c(\Z)>0$ and $E(\Z, \Z)=0$, then $\tors(\Z)$ is an infinite set.
\end{enumerate}

\end{proposition}

\begin{proof}
    Part $(1)$ follows from the proof of the forward direction of Proposition \ref{Prop: infinite tors(Z)}. 
    For part $(2)$, assume that $\Z$ is $\tau$-regular with $c(\Z)>0$ and $E(\Z, \Z)=0$.

     If $(Z_1, Z_2) \in \Z \times \Z$ is in general position, then $E(\Z,\Z)=0$ means that $\Hom(Z_1, \tau Z_2)=0$. Since $c(\Z) \ge 1$, the $Z_i$ can be taken so that $\Hom(Z_i, \tau Z_i) = c(\Z) \ne 0$. We let $B_i$ be the image of a non-zero morphism $Z_i \to \tau Z_i$ which has minimal dimension. We let $\pi_i: Z_i \to B_i$ be the canonical projection and $\iota_i: B_i \to \tau Z_i$ be the canonical injection. If $\varphi_i: B_i \to B_i$ is nonzero and not-invertible, the map $\iota_i\varphi_i\pi_i: Z_i \to \tau Z_i$ is non-zero with image a proper submodule of $B_i$, which contradicts the minimality of $B_i$. Hence, $B_i$ is a brick. Similarly, if $\varphi_{i,j}: B_i \to B_j$ is non-zero, then $\iota_j\varphi_{i,j}\pi_i: Z_i \to \tau Z_j$ is non-zero, a contradiction. This shows that $B_i, B_j$ are Hom-orthogonal bricks. The fact that $c(\Z) \ge 1$ allows us to start with an infinite family $\{Z_\lambda\}_{\lambda \in \Lambda}$ of modules in general position in $\Z$ such that $\Hom(Z_\mu, \tau Z_{\nu})=0$ for $\mu \ne \nu$. Applying the above argument gives rise to an infinite family $\{B_\lambda\}_{\lambda \in \Lambda'}$ of Hom-orthogonal bricks. 
\end{proof}

\begin{remark} \label{Rem: stability}
    In the above proof, we constructed bricks $\{B_\lambda\}$ from the $\{Z_\lambda\}$, the latter lying in the $\tau$-regular component $\Z$. Consider  $Z_\lambda$ and $Z_\mu$ with $\mu, \nu \in \Lambda$, and consider the corresponding bricks $B_\lambda, B\mu$. Let $v$ be the $g$-vector corresponding to $\Z$. We can realize $v$ in two different ways:
    \[v = {\rm dim}_k \Hom_A(Z_\lambda, -) - {\rm dim}_k\Hom_A(-, \tau Z_\lambda).\]
    \[v = {\rm dim}_k \Hom_A(Z_\mu, -) - {\rm dim}_k\Hom_A(-, \tau Z_\mu).\]
    We have $\Hom(Z_\mu, B_\lambda)=0$ and $\Hom(B_\lambda, \tau Z_\mu)=0$. Indeed, this follows from the fact that $\Hom(Z_\mu, \tau Z_\lambda)= \Hom(Z_\lambda, \tau Z_\mu)=0$. Using the second interpretation of $v$, this yields $v([B_\lambda])=0$. Now, if $L$ is any proper submodule of $B_\lambda$, then $\Hom(Z_\lambda, L)=0$, by the minimality property in the definition of $B_\lambda$. However, the canonical inclusions $L \hookrightarrow B_\lambda \hookrightarrow \tau Z_\lambda$ imply that $\Hom(L, \tau Z_\lambda) \ne 0$. This shows, using the first interpretation of $v$, that $v([L])<0$. Hence, all $B_\lambda$ are in fact $v$-stable.
\end{remark}

\subsection{$E$-tame algebras}

We now turn our attention to $E$-tame algebras. Recall that $A$ is \emph{$E$-tame} if for every indecomposable $g$-vector $\theta$, and each positive integer $r$, the canonical decomposition of $r\theta$ is 
$$r\theta = \theta + \cdots + \theta.$$
Such algebras manifest many important properties and they have received a lot of attention in recent years (For details, see \cite{AI}, \cite{PY}, and references therein). Let us briefly remark that the family of $E$-tame algebras strictly generalizes the family of tame algebras; in particular, every brick-finite algebra is $E$-tame. 

\medskip

The following result directly follows from Proposition \ref{Prop on c(Z) zero}, hence we omit the proof. Observe that the next statement essentially says that, for $E$-tame algebras, Demonet's Conjecture implies the 2nd bBT Conjecture. Recall that these two conjectures are known to be equivalent over tame algebras (see Remark \ref{2nd bBT and Demonet's conj for tame algebras}).

\begin{corollary} \label{Cor: E-tame and bricks}
    Let $A$ be an $E$-tame algebra that has a $\tau$-regular component $\Z$ with $c(\Z)>0$ (equivalently, there is a integral $g$-vector outside the $\tau$-tilting fan). Then, $A$ admits a brick component $\Z'$ with $c(\Z')>0$.
\end{corollary}

We end this section with the following proposition, which is a more general version of Corollary \ref{Cor:tame-almost all bricks tau-rigid}. We again remark that any tame algebra is $E$-tame. 

\begin{proposition}\label{prop: E-tame almost all brick tau-rigid}
    Let $A$ be $E$-tame such that almost all bricks are $\tau$-rigid. Then $A$ is brick-finite.
\end{proposition}

\begin{proof}
The proof is the same as that of Proposition \ref{Cor:tame-almost all bricks tau-rigid}, knowing that, thanks to Corollary \ref{Cor: E-tame and bricks}, for an $E$-tame algebra, having a rational vector outside $\mathcal{F}_A$ implies the existence of an infinite family of bricks in the same dimension.
\end{proof}

\textbf{Acknowledgements.} 
KM was supported by Early-Career
Scientist JSPS Kakenhi grant number 24K16908. CP was supported by the Natural Sciences and Engineering Research Council of Canada and by the Canadian Defence
Academy Research Programme. The earlier part of this project was conducted while the first-named author visited Bielefeld Representation Theory group (BIREP) for $6$ weeks, facilitated by the SFB/TRR 358 funding. He is grateful to the BIREP group for their hospitality, and particularly thanks William Crawley-Boevey and Henning Krause for several stimulating discussions. The authors also thank Sota Asai, Raymundo Bautista, Christof Geiß, Osamu Iyama, Rosanna Laking and Claus Michael Ringel, for helpful communications at different stages of this project.

\end{document}